\documentclass[12pt]{amsart}
\usepackage[a4paper,margin=2cm]{geometry}

\usepackage[T1]{fontenc} 
\usepackage[utf8]{inputenc}

\usepackage{amsmath,amssymb,amsthm,mathtools,enumerate,extarrows,stmaryrd,mathrsfs,textcomp}
\usepackage[oldstyle]{libertine}

\usepackage[dvipsnames]{xcolor}

\usepackage{enumitem}
\newcommand{\desclabel}[2]{%
  \item[#1]\label{#2}\def\@currentlabelname{#1}%
}

\usepackage{graphicx} 
\usepackage{subfig}

\newcommand{\descref}[1]{\hyperref[#1]{\nameref{#1}}}

\usepackage[dvipsnames]{xcolor}

\usepackage{mathtools}
\usepackage{extarrows}
\usepackage{tikz}
\usepackage[bbgreekl]{mathbbol}
\DeclareSymbolFontAlphabet{\mathbbm}{bbold}
\DeclareSymbolFontAlphabet{\mathbb}{AMSb}%
   
\DeclareMathOperator*{\Res}{Res}

\UseRawInputEncoding
\let\oldtocsubsection=\tocsubsection
\renewcommand{\tocsubsection}[2]{\hspace{1em}\oldtocsubsection{#1}{#2}}

\usepackage{hyperref}
\hypersetup{colorlinks=true,linkcolor=blue,citecolor=teal,filecolor=magenta,urlcolor=cyan}

\usepackage{url}

\usepackage[backend=bibtex,style=alphabetic,sorting=nyt,isbn=false,url=false,doi=false,maxalphanames=10,minalphanames=4,mincitenames=4,maxcitenames=10,minnames=4,maxnames=10,giveninits=true,maxbibnames=99]{biblatex}
\usepackage[newparttoc]{titlesec}
\usepackage{appendix}
\usepackage{titletoc}
\titleformat{\section}
{\normalfont\bfseries\large}{\thesection}{1em}{}
\titlecontents{section}
	[3.5em]{}
	{\contentslabel{1.5em}}
	{}
	{\hfill\contentspage}
\titleformat{\subsection}
{\normalfont\large\slshape}{\thesubsection}{1em}{}
\titlecontents{subsection}
	[5.5em]{}
	{\contentslabel{1.5em}\itshape}
	{}
	{\hfill\contentspage}
\titleformat{\subsubsection}
{\normalfont\itshape}{\thesubsubsection}{1em}{}
\titlecontents{subsubsection}
	[7.5em]{}
	{\contentslabel{1.5em}\itshape}
	{}
	{\hfill\contentspage}

\theoremstyle{plain}
\newtheorem{theorem}{Theorem}[section]
\newtheorem{proposition}[theorem]{Proposition}

\newtheorem{question}[theorem]{Question}
\newtheorem{conjecture}[theorem]{Conjecture}
\newtheorem{lemma}[theorem]{Lemma}
\theoremstyle{definition}
\newtheorem{definition}[theorem]{Definition}

\theoremstyle{remark}
\newtheorem{remark}[theorem]{Remark}

\numberwithin{equation}{section}

\makeatletter
\AtBeginEnvironment{abstract}{%
  \setbox0=\hbox{$0$}
  \textfont9=\textfont1
  \scriptfont9=\scriptfont1
  \scriptscriptfont9=\scriptscriptfont1
}
\makeatother

\begin{document}
	
\title{Volumes of moduli spaces of bordered Klein surfaces}

\author[E.~Garcia-Failde]{Elba Garcia-Failde}
\address[E.~G.]{Institut de Math\'ematiques de Jussieu - Paris Rive Gauche \\ Sorbonne Universit\'e, Paris, France.}

\address[E.~G.]{Departament de Matem\`{a}tiques  \\ Universitat Polit\`{e}cnica de Catalunya, Barcelona, Spain.}
\email[E.~G.]{\href{mailto:elba.garcia@upc.edu}{elba.garcia@upc.edu}}

\author[P.~Gregori]{Paolo Gregori}
\address[P.~G.]{Universit\'e Paris-Saclay, CNRS, CEA, Institut de Physique Th\'eorique, Gif-sur-Yvette, France.}
\email[P.~G.]{\href{mailto:plgregori@gmail.com}{plgregori@gmail.com}}

\author[K.~Osuga]{Kento Osuga}
\address[K.~O.]{Graduate School of Mathematical Sciences, University of Tokyo, Tokyo, Japan.}
\address[K.~O.]{Kobayashi--Masakawa Institute for the Origin of Particles and the Universe}
\address[K.~O.]{Graduate School of Mathematics, Nagoya University, Nagoya, Japan.}
\email[K.~O.]{\href{mailto:osuga@math.nagoya-u.ac.jp}{osuga@math.nagoya-u.ac.jp}}

\begin{abstract}


We investigate volumes of moduli spaces of bordered Klein surfaces, which include non-orientable surfaces. 
On these moduli spaces, the top form introduced by Norbury diverges as the lengths of 1-sided geodesics approach zero. However, when integrated over Gendulphe's regularised moduli space, on which the systole of 1-sided geodesics is bounded below by $\epsilon\in\mathbb{R}_{>0}$, it returns a finite value.
We derive an explicit formula for the volume of the moduli space of two-bordered real projective planes and of one-bordered Klein bottles --- the latter utilises Norbury's extension of the Mirzakhani--McShane identities to non-orientable surfaces. We then relate these results to refined topological recursion, showing that, for a fixed refinement parameter, the volumes of moduli spaces of Klein surfaces with Euler characteristic $-1$ are governed by this procedure. For general topologies, both approaches remain incomplete due to several difficulties, which we discuss in this work. Together with our results, this initiates the search for a geometric recursive structure governing the volumes of moduli spaces of bordered Klein surfaces, whose natural formulation may lie within a suitably extended version of refined topological recursion.

\end{abstract}

\maketitle

\tableofcontents	

\newpage

\section{Introduction}\label{sec:intro}

The aim of this article is twofold: one is to compute volumes of moduli spaces of bordered Klein surfaces of Euler characteristic $-1$, and the other is to explore their recursive structure for general Euler characteristics from the perspective of refined topological recursion. The former is of interest in hyperbolic geometry and topology, whereas the latter plays a role in the context of mirror symmetry. These two perspectives also lead us to an intriguing one-parameter deformation that continuously interpolates between orientable and non-orientable surfaces.

\subsection{Motivation and background}

For $g\in\mathbb{Z}_{\geq0}$ and $n\in\mathbb{Z}_{>0}$ with $2-2g-n<0$, let $\mathcal{M}^+_{g,n}(L_{[n]})$ be the moduli space of oriented hyperbolic surfaces of genus $g$ with $n$ boundaries of lengths $L_{[n]}=(L_1,...,L_n)$. It is well-known that $\mathcal{M}^+_{g,n}(L_{[n]})$ admits a symplectic form $\omega_{{\rm WP}}$, called the Weil--Petersson form, from which one may define its symplectic volume $V^+_{g,n}(L_{[n]})$, called the \emph{Weil--Petersson volume}. Mirzakhani proved in \cite{Mir06} that the Weil--Petersson volumes $V^+_{g,n}(L_{[n]})$ satisfy a recursive formula, now known as \emph{Mirzakhani's recursion}. It is a recursion on the negative Euler characteristic of the surface, $2g-2+n$, and admits an interpretation in terms of pants decompositions. Mirzakhani derived such a recursion by applying the so-called Mirzakhani--McShane identities \cite{McS98,Mir06}, which classify the fates of ortho-geodesics emanating from the first boundary component for every fixed $(g,n)$. It was further shown in \cite{Mir07,MS08} that $V^+_{g,n}(L_{[n]})$ is related to $\psi$-class intersection numbers on the moduli space of stable curves, Virasoro constraints, and a $\tau$-function of the Korteweg--de Vries hierarchy.  

Shortly after, Chekhov, Eynard, and Orantin showed that Mirzakhani's recursion could be reformulated in terms of a universal procedure called \emph{topological recursion} \cite{CE05,EO07}.
The set of initial data forms a spectral curve $\mathcal{S}_\Sigma$, which consists of a Riemann surface $\Sigma$ together with a choice of unstable correlators $(\omega_{0,1},\omega_{0,2})$. The output is a sequence of stable correlators $(\omega_{g,n})_{2g-2+n>0}$, which are meromorphic sections of $K_\Sigma^{\boxtimes n}$. Although it was originally discovered in the context of Hermitian matrix integrals, it was shown in \cite{EO08} that Mirzakhani's recursion is actually an instance of topological recursion on a specific spectral curve $\mathcal{S}_{\rm WP}$.  It turns out that topological recursion can be applied to far more general enumerative problems, e.g.~Gromov--Witten invariants of a toric Calabi-Yau threefold $X$ \cite{BKMP08,FLZ16}. In the context of mirror symmetry, topological recursion is thus considered as the B-side mirror of Gromov--Witten theory on $X$, since the corresponding spectral curve $\mathcal{S}_X$ is the Hori--Vafa mirror curve of $X$ \cite{HV00}.

Recently, a one-parameter deformation of topological recursion, called \emph{refined topological recursion}, was proposed in \cite{KO22,Osu23-1,Osu23-2} inspired by $\beta$-deformed matrix integrals \cite{CE06}. The set of initial data is upgraded to a refined spectral curve $\mathcal{S}_\Sigma^\mathfrak{b}$, which includes unstable correlators $\omega_{g,n}^\mathfrak{b}$ with $(g,n)=(0,1), (0,2)$, and additionally $(\frac12,1)$ carrying a deformation parameter $\mathfrak{b}$. The output is a sequence of refined stable correlators $(\omega^\mathfrak{b}_{g,n})_{2g-2+n>0}$, labelled by an integer $n\in\mathbb{Z}_{>0}$ and a half-integer $g\in\frac12\mathbb{Z}_{\geq0}$. It has then been shown in \cite{CDO24-1,CDO24-2} that refined correlators $\omega^\mathfrak{b}_{g,n}$ on a certain refined spectral curve $\mathcal{S}^\mathfrak{b}_{\rm Hur}$ are related to $b$-Hurwitz numbers $H_{g,n}^b$, which count embedded graphs $\Gamma$ on possibly non-orientable surfaces, weighted by a \emph{measure of non-orientability} $\rho_\Gamma(b)$ introduced by Chapuy and Do\l\k{e}ga \cite{CD20}. The parameters $b$ and $\mathfrak{b}$ are related by $\mathfrak{b}=-\frac{b}{\sqrt{1+b}}$, and $b$ intuitively measures ``how non-orientable'' the embedded graph $\Gamma$ is, in a suitable sense --- in particular, when $b=0$, only orientable contributions remain nonzero.

\medskip

Motivated by these recent developments, intriguing questions arise:

\begin{itemize}
\item[\textbf{Q1:}] Can one consider volumes of moduli spaces of non-orientable surfaces, and extend Mirzakhani's recursion to the non-orientable setting?
\item[\textbf{Q2:}] Can such a non-orientable analogue of Mirzakhani's recursion be deformed by means of the parameter $b$, and can this deformation be interpreted geometrically in terms of a measure of non-orientability?
\end{itemize}
In this article, we report partial progress on both \textbf{Q1} and \textbf{Q2}, which we summarise below.

\subsection{Main results and future directions}

In \cite{Nor07}, Norbury introduced a top-degree form $\nu_{\rm N}$, which we call the Norbury form, on the moduli space $\mathcal{M}^-_{g,n}(L_{[n]})$ of non-orientable hyperbolic surfaces of Euler characteristic $2-2g-n$, with $n$ boundaries of lengths $L_{[n]}$, and where $g\in\frac12\mathbb{Z}_{>0}$ denotes \emph{half} the non-orientable genus (some technical details are omitted for now). Since the Norbury form $\nu_{\rm N}$ is invariant under the action of the mapping class group, one may wish to take it as a non-orientable analogue of the Weil--Petersson volume form. However, it turns out that $\nu_{\rm N}$ has singularities on $\mathcal{M}^-_{g,n}(L_{[n]})$, hence one faces the immediate issue that the naive integral of the Norbury form $\nu_{\rm N}$ diverges. 

Gendulphe \cite{Gen17} resolved this issue by introducing a \emph{regularised moduli space} $\mathcal{M}^{-,\epsilon}_{g,n}(L_{[n]})$, on which the Norbury form remains finite for a sufficiently small regularisation parameter $\epsilon\in\mathbb{R}_{>0}$. Gendulphe's regularisation is geometrically natural: $\epsilon$ is interpreted as the lower bound of the systole of 1-sided curves, regardless of the topology. Therefore, one can now define $V_{g,n}^{-,\epsilon}(L_{[n]})$, the non-orientable analogue of the Weil--Petersson volume, by the integral of the Norbury form $\nu_{\rm N}$ over Gendulphe's regularised moduli space $\mathcal{M}^{-,\epsilon}_{g,n}(L_{[n]})$, with a suitable normalisation constant. We call $V_{g,n}^{-,\epsilon}(L_{[n]})$ the \emph{Gendulphe--Norbury volume}.

As a trade-off for being able to define finite volumes $V_{g,n}^{-,\epsilon}$, the regularised moduli spaces $\mathcal{M}^{-,\epsilon}_{g,n}(L_{[n]})$ become more complicated, which makes concrete computations difficult. This is partly why no significant progress has been made after the introduction of Gendulphe's regularisation. In order to break this impasse, we take the first step towards the full understanding of $V_{g,n}^{-,\epsilon}$:
\begin{theorem} The Gendulphe--Norbury volumes of the moduli spaces of two-bordered real projective planes $(g,n)=(\frac12,2)$ and of one-bordered Klein bottles $(g,n)=(1,1)$ are given as:
\begin{align*}
V_{\frac12,2}^{-,\epsilon}(L_{[2]})&=\log\frac{\cosh\frac{L_1+L_2}{4}\cosh\frac{L_1-L_2}{4}}{\sinh^2\frac{\epsilon}{2}},\\
V_{1,1}^{-,\epsilon}(L_1)&=-{\rm Li}_2\left(-\frac{\cosh^2\frac{L_1}{4}}{\sinh^2\frac{\epsilon}{2}}\right).
\end{align*}
\end{theorem}

Note that both of them blow up as $\epsilon\to0$ as expected. Let us emphasise that, since $\mathcal{M}^{-,\epsilon}_{g,n}(L_{[n]})$ are complicated due to Gendulphe's regularisation, it is remarkable in particular that the volume for one-bordered Klein bottles $V_{1,1}^{-,\epsilon}(L_1)$ admits such a simple expression in terms of the dilogarithm function. One should also notice a clear contrast to the orientable case in which all the Weil--Petersson volumes $V_{g,n}^{+}(L_{[n]})$ are polynomials in $L_{[n]}$, while polynomiality is drastically broken in the non-orientable setting. 
While the computation for $(g,n)=(\frac12,2)$ has been previously obtained in the physics context (see~\cite{Sta23}), our derivation offers a geometric perspective, and the result for $(g,n)=(1,1)$ is new.

\medskip

We will now change our perspective to refined topological recursion. Since the formalism itself was already proposed in \cite{KO22,Osu23-1}, a key challenging task is to determine the spectral curve $\mathcal{S}^\mathfrak{b}_{\rm total}$ whose associated refined correlators $\omega_{g,n}^\mathfrak{b}$ recover the \emph{total volumes} \[V_{g,n}^{\epsilon}(L_{[n]}):=V_{g,n}^{+}(L_{[n]})+V_{g,n}^{-,\epsilon}(L_{[n]})\,,\]
where we conventionally set $V_{g,n}^{+}(L_{[n]})=0$ whenever $g\in\mathbb{Z}_{\geq0}+\frac12$. In addition to the data included in $\mathcal{S}_{\rm WP}$, which in particular sets $\Sigma=\mathbb{C}$, one needs to find $\omega^\mathfrak{b}_{g,n}$ for $(g,n)=(\frac12,1)$. Guided by extensive computational experiments, we found the following well-supported candidate:
\[\omega_{\frac12,1}^\mathfrak{b}(z)\coloneqq\frac{\mathfrak{b}}{2}\left(-\frac{2\pi}{\tan2\pi z}+\frac{e^{2z\epsilon}-e^{-2z\epsilon}}{2z}+\sum_{k\geq1}\left(\frac{e^{2(z-\frac{k}{2})\epsilon}}{z-\frac{k}{2}}-\frac{e^{-2(z+\frac{k}{2})\epsilon}}{z+\frac{k}{2}}\right)\right)dz,\]
where $z$ is the standard coordinate on $\mathbb{C}$. The summation is convergent when $\epsilon\in\mathbb{R}_{>0}$, matching one of the conditions required for Gendulphe's regularisation. 

Our second main result proves that refined stable correlators $\omega_{g,n}^\mathfrak{b}$ constructed on $\mathcal{S}^\mathfrak{b}_{\rm total}$ recover total volumes $V_{g,n}^{\epsilon}$ when $(g,n)=(0,3), (\frac12,2), (1,1)$. Omitting some details, we propose the following concrete dictionary between the refined correlators $\omega_{g,n}^\mathfrak{b}$ and the total volumes $V_{g,n}^{\epsilon}$:

\begin{theorem}\label{thm:intro3}
For $(g,n)\in\frac12\mathbb{Z}_{\geq0}\times \mathbb{Z}_{\geq1}$ with $2-2g-n=-1$, the refined correlators $\omega_{g,n}^\mathfrak{b}$ and the total volumes $V_{g,n}^{\epsilon}$ are related, in a subset of $\mathbb{R}_{>0}^n$, by
\[  \prod_{i=1}^nL_i\cdot V_{g,n}^\epsilon(L_{[n]})= (1+b)^{g}\cdot (\hat{\mathcal{L}}_{L_1}^{-1}\cdots \hat{\mathcal{L}}_{L_n}^{-1}).\;\omega^\mathfrak{b}_{g,n}\Big|_{b=1},\]
where $\hat{\mathcal{L}}^{-1}$ is the termwise inverse Laplace transform and $\mathfrak{b}=-\frac{b}{\sqrt{1+b}}$.
\end{theorem}

The overall factor $(1+b)^g$ is inspired by the relation between refined topological recursion and $b$-Hurwitz numbers \cite{CDO24-1,CDO24-2}. We note that establishing this correspondence is highly nontrivial, particularly in the case $(g,n)=(1,1)$, where the two sides initially take strikingly different forms and are shown to coincide only after suitable expansions.

\subsubsection{Future direction 1: general topologies}

So what about topologies with $2-2g-n<-1$? From the perspective of refined topological recursion, problems originate from the fact that all properties proved in \cite{KO22,Osu23-1} hold when $\Sigma=\mathbb{P}^1$, but not when $\Sigma=\mathbb{C}$. For example, it is not clear whether the refined recursion formula generates symmetric $\omega_{g,n}^\mathfrak{b}$. Furthermore, even in this case, they are no longer rational differentials. Instead, they are, in general, divergent series, which are expected to become convergent after a termwise inverse Laplace transform. These features are verified  in explicitly provable cases, i.e.~$(g,n)=(0,3),(\frac12,2),(1,1)$, but remain open for the other $\omega_{g,n}^\mathfrak{b}$ at the time of writing. 

In hyperbolic geometry, on the other hand, Norbury already derived non-orientable analogues of McShane--Mirzakhani identities for general topologies \cite{Nor07}, which we call Norbury identities. Indeed, we used the Norbury identity for one-bordered Klein bottles to derive the expression for $V_{1,1}^{-,\epsilon}(L_1)$ above. For $2-2g-n<-1$, one may hope to apply the Norbury identity for the corresponding topology to derive a recursive equation for the Gendulphe--Norbury volumes $V_{g,n}^{-,\epsilon}(L_{[n]})$, in close analogy with Mirzakhani's approach \cite{Mir06}. However, one notices that Gendulphe's regularisation is difficult to keep track of in a recursive procedure --- new 1-sided curves will be created by gluing a pair of pants to a Klein surface, say $K$, and their lengths can be smaller than $\epsilon$ even if all 1-sided curves in $K$ are bounded below by $\epsilon$\footnote{We acknowledge Bram Petri for raising questions about this crucial feature, following discussions with A.~Wright.}.

Despite these issues on both sides, the available evidence suggests that Theorem~\ref{thm:intro3} may extend to general topologies. Indeed, our analysis indicates that a clean recursion involving only the total volumes $V_{g,n}^{\epsilon}(L_{[n]})$ is unlikely to exist, whereas a more tractable recursive formula may emerge when one instead takes the refined correlators $\omega_{g,n}^\mathfrak{b}$ as the fundamental objects. It is also worth emphasising that, whereas the quantities $V_{g,n}^{\epsilon}(L_{[n]})$ with $2g-2+n=1$ would have to be supplied as initial data for such a recursion (were one to exist), the refined correlators $\omega_{g,n}^{\mathfrak{b}}$ with $2g-2+n=1$ are themselves determined by refined topological recursion, since its input consists only of the cases with $2g-2+n\leq0$. Thus, refined topological recursion has already passed at least one nontrivial check for the existence of such a recursive formula.

Furthermore, \cite{CGGO26} establishes a relation between refined topological recursion and the volumes $V^{\text{M\"{o}}}_{g,n}(L_{[n]})$ of moduli spaces of \emph{metric M\"{o}bius graphs}, which is expected to be equivalent to Theorem~\ref{thm:intro3} in the limit $L_{k}\to\infty$. Unlike the hyperbolic geometry setting, both sides of this correspondence are shown to admit recursive formulas in \cite{CGGO26}. However, the recursion for $V^{\text{M\"{o}}}_{g,n}(L_{[n]})$ is considerably less transparent because these volumes are piecewise polynomial, a feature that may be viewed as a remnant of the intricate structure of $V_{g,n}^{-,\epsilon}(L_{[n]})$. In contrast, the recursion for $\omega_{g,n}^\mathfrak{b}$ remains simpler and universal, independently of the moduli spaces under consideration. 

Although this question lies beyond the scope of the present paper, we hope to return to it in the near future.

\subsubsection{Future direction 2: measure of non-orientability}

Having summarised a (conjectural) relation between the total volumes and refined topological recursion, a natural question arises: what happens if the parameter $b$ is left unfixed in Theorem~\ref{thm:intro3}? To be more precise, let us \emph{define} refined total volumes $V_{g,n}^{\epsilon,b}$ by
\[  \prod_{i=1}^nL_i\cdot V_{g,n}^{\epsilon,b}(L_{[n]})\coloneqq (1+b)^{g}\cdot (\hat{\mathcal{L}}_{L_1}^{-1}\cdots \hat{\mathcal{L}}_{L_n}^{-1}).\;\omega^\mathfrak{b}_{g,n},\]
for all $(g,n)$ with $2-2g-n<0$. It is easy to see that, under suitable technical assumptions, $V_{g,n}^{\epsilon,b}$ depends polynomially on $b$. This naturally leads to the following question: can we assign a geometric meaning to $V_{g,n}^{\epsilon,b}$ by means of a measure of non-orientability analogous to the one introduced by Chapuy--Do\l\k{e}ga for $b$-Hurwitz numbers \cite{CD20}? For instance, do the coefficients of $b$ in $V_{1,1}^{\epsilon,b}$ below admit geometric interpretations? This question deserves further investigation.
\[V_{1,1}^{\epsilon,b}=\frac{L^2}{48}+\frac{\pi^2}{12}+b\left(\frac{L^2}{48}+\frac{\pi^2}{12}\right)-b^2\left(\frac{L^2}{48}+\frac{\pi^2}{12}+{\rm Li}_2\left(-\frac{\cosh^2\frac{L_1}{4}}{\sinh^2\frac{\epsilon}{2}}\right)\right).\]


\subsection*{Organisation}

This paper is organised as follows. Section~\ref{sec:non-orientable geometry} reviews hyperbolic geometry on non-orientable surfaces and introduces the Norbury form and Gendulphe's regularisation. We also recall Mirzakhani's unfolding technique, which can be partially extended to the non-orientable setting. Section~\ref{sec:euler characteristic -1 from geometry} then specialises to surfaces of Euler characteristic $2-2g-n=-1$, where we explicitly compute $V_{g,n}^{-,\epsilon}$ for $(g,n)=(\frac12,2)$ and $(1,1)$, establishing one of the main results of the paper.
Section~\ref{sec:higherTop} discusses obstructions to a straightforward extension of this geometric approach to general topologies, motivating the introduction of refined topological recursion. We then investigate its relation to the total volumes, which may be viewed as a manifestation of mirror symmetry in the refined setting. Section~\ref{sec:evidence} provides evidence for this mirror correspondence by establishing that $\omega_{g,n}^{\mathfrak b}$ and $V_{g,n}^{\epsilon}(L_{[n]})$ are related through an explicit dictionary in the cases with $2g-2+n=1$, constituting the second main result of the paper. Appendix~\ref{sec:formulae} collects several combinatorial formulae used to establish the correspondence in the case $(g,n)=(1,1)$. Finally, Appendix~\ref{sec:subtleties_in_question_ref_conj_main} discusses a subtlety that arises in a naive attempt to relate refined topological recursion to the total volumes.


\subsection*{Acknowledgements}

The authors would like to thank M.~Gendulphe, P.~Norbury and B.~Petri for helpful comments on this work. We also thank N.~Chidambaram, A.~Giacchetto and D.~Lewa\'nski for valuable discussions. We are particularly grateful to B.~Petri and A.~Wright for raising fundamental questions about earlier versions of the recursion. Finally, we thank the Institut de Physique Th\'eorique in Saclay, the Institut de Math\'ematiques de Jussieu - Paris Rive Gauche in Sorbonne Universit\'e, the Research Institute for Mathematical Sciences at Kyoto University, and the Graduate School of Mathematics at the University of Tokyo for their hospitality.

E.~G-F.~is currently supported by a Ram\'on y Cajal fellowship. She also acknowledges the support of a Tremplin grant from Sorbonne Universit\'e, a PEPS grant from the CNRS, the ERC Synergy Grant ReNewQuantum and the project PID2024-155686NB-I00 of the Spanish Ministry of Science and Innovation.

P. G.~was supported by the ERC Synergy Grant ReNewQuantum.

K.~O.~acknowledges the support by JSPS KAKENHI Grant-in-Aid for JSPS Fellows (22KJ0715) and  for Early-Career Scientists (23K12968, 26K16980), and in part also by 24K00525. K.~O.~also acknowledges the support from the Kobayashi--Maskawa Institute (KMI) for the Origin of Particles and the Universe at Nagoya University.  


\newpage 

\section{Non-orientable hyperbolic geometry}\label{sec:non-orientable geometry}

We review the necessary background on hyperbolic geometry for Klein surfaces. We refer the reader to \cite{Jos97,Nor07,Bus10,PP13,Gen17} for further details.

We first summarise some topological aspects. For an integer $n\in\mathbb{Z}_{\geq1}$ and a half-integer $g\in\frac12\mathbb{Z}_{\geq0}$, let $K$ be an $n$-bordered surface of genus $g$\footnote{We abuse notation for simplicity throughout the paper; $g$ is the genus if $K$ is orientable, but it is \emph{half} of the non-orientable genus if $K$ is non-orientable. With this notation, the Euler characteristic reads $\chi(K)=2-2g-n$ regardless of the orientability of $K$.} which may or may not be orientable. A simple curve $\gamma\subset K$ is called \emph{1-sided} (resp.~\emph{2-sided}) if its tubular neighbourhood is a M\"{o}bius strip (cylinder). A 2-sided simple curve is called \emph{non-primitive} if and only if it bounds a M\"{o}bius strip --- whose core is of course 1-sided ---, and \emph{primitive} otherwise. For brevity of terminology, all curves in the present paper are simple, and 2-sided curves always mean primitive unless specified. A \emph{pants decomposition} $\mathcal{P}\subset K$ is a set of non-boundary parallel 1-sided or 2-sided curves, up to homotopy, such that $K\backslash \mathcal{P}$ is a disjoint union of pairs of pants.

For an index set $I$, we denote by $X_{I}=(X_i)_{i\in I}$ the family of elements $X_i$. Then, for $N_1,N_2\in\mathbb{Z}_{\geq0}$ with $N_1+2N_2 = 6g-6+2n\geq0$, one can decompose a pants decomposition $\mathcal{P}$ into a set of 1-sided curves $\alpha_{[N_1]}$ and 2-sided curves $\gamma_{[N_2]}$, where $[N_1]=\{1,...,N_1\}$ and $[N_2]=\{1,...,N_2\}$, with $[0]=\varnothing$. Note that a pants decomposition $\mathcal{P}=\alpha_{[N_1]}\sqcup \gamma_{[N_2]}$ is not unique, in particular, $N_1$ and $N_2$ are not topologically invariant numbers, and $N_1$ is bounded above by $2g$. Indeed, up to homotopy, two pants decompositions are related by a composition of the four elementary moves below \cite{PP13}:
\begin{enumerate}
\item Replace a 2-sided curve $\gamma$ in a four-bordered sphere with another 2-sided curve $\gamma'$ intersecting $\gamma$ exactly once.
\item Replace a 2-sided curve $\gamma$ in a one-bordered torus with another 2-sided curve $\gamma'$ intersecting $\gamma$ exactly once.
\item Replace a 1-sided curve $\alpha$ in a two-bordered real-projective plane with another 1-sided curve $\alpha'$ intersecting $\alpha$ exactly once.
\item Replace two non-intersecting 1-sided curves $\alpha,\alpha'$ in a one-bordered Klein bottle with a two-sided curve $\gamma$ intersecting both $\alpha,\alpha'$ exactly once each.
\end{enumerate}

\begin{figure}[ht]
  \centering
  \subfloat{\includegraphics[width=0.3\textwidth]{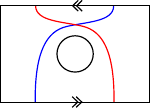}}
  \hskip 50pt  
  \subfloat{\includegraphics[width=0.324\textwidth]{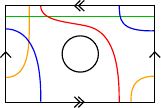}}
  \caption{The picture on the left shows a two-bordered real-projective plane which contains exactly two 1-sided curves $\alpha$ and $\alpha'$ (red and blue) intersecting exactly once; the third elementary move exchanges them. The picture on the right shows a one-bordered Klein bottle with its unique 2-sided primitive curve (green), a reference 1-sided curve $\alpha_0=\alpha$ (red), and a pair of 1-sided curves $\alpha_1=\alpha'$ and $\alpha_{-1}$ (blue and orange) which intersect each other once but do not intersect $\alpha_0$, and can be obtained respectively by adding/subtracting a Dehn twist along $\gamma$ to $\alpha_0$. The fourth elementary move exchanges the pair $(\alpha,\alpha')$ with $\gamma$, which intersects both exactly once.}\label{fig:moves}
\end{figure}


\subsection{Gendulphe--Norbury's regularised volumes}\label{vols}

We next endow $K$ with a hyperbolic metric $\rho$ and set the boundaries to be geodesics of lengths $L_{[n]}$. A pair $(K,\rho)$ is often called a \emph{Klein surface} in the literature, when $L_i=0$ for all $i\in[n]$\footnote{A Klein surface can be realised as a dianalytic manifold of complex dimension 1. We recall that a \emph{dianalytic structure} is given by an atlas whose changes of charts are holomorphic or anti-holomorphic.}, and when $L_i>0$, we call a pair $(K,\rho)$ a \emph{Klein surface with boundaries}. When $K$ is orientable and oriented, it is simply a Riemann surface with $n$ boundaries. The \emph{Teichm\"{u}ller space} $\mathcal{T}(K)$ is the space of isotopy classes of hyperbolic metrics $\rho$ on $K$. It is known that $\dim_\mathbb{R}\mathcal{T}(K)=6g-6+2n$.

Given a pants decomposition $\mathcal{P}\subset K$, we consider for each curve in $\mathcal{P}$ the shortest-geodesic representative in its homotopy class. We denote their lengths by $\ell_{\alpha_1},...,\ell_{\alpha_{N_1}}$ for 1-sided curves and $\ell_{\gamma_1},...,\ell_{\gamma_{N_2}}$ for 2-sided curves. It is known (e.g.~\cite[Theorem 2.2]{PP13}) that $\mathcal{T}(K)$ is homeomorphic (and, in fact, real-analytically diffeomorphic) to $(\mathbb{R}_{>0})^{N_1}\times(\mathbb{R}_{>0}\times\mathbb{R})^{N_2}$, whose coordinates are given by $\ell_{\alpha_1},...,\ell_{\alpha_{N_1}}$ and $\ell_{\gamma_1},...,\ell_{\gamma_{N_2}}$ together with \emph{twisting coordinates} $\theta_1,...,\theta_{N_2}$. We call $(\ell_{\alpha_i},\ell_{\gamma_j},\theta_j)^{i\in[N_1]}_{j\in[N_2]}$ the \emph{Fenchel--Nielsen coordinates} (for Klein surfaces with boundaries).

Twisting coordinates represent the angle upon gluing of two 2-sided curves. However, to make sense of the twisting parameter $\theta_j$ at each 2-sided geodesic $\gamma_j$, one needs a certain notion of orientation. Let us now clarify this, closely following \cite{Nor07}. 

\subsubsection{Twisting coordinates and orientation of the first boundary component}\label{sec:twist coordinates}

When $K$ is orientable, we simply pick an orientation on it and $\mathcal{T}(K)$ refers to the Teichm\"{u}ller space of oriented hyperbolic surfaces. In this case, when cutting $K$ along the curves in $\mathcal{P}=\gamma_{[N_2]}$, each boundary in the pairs of pants inherits an orientation from the orientation of $K$ and one can reconstruct the original surface by pairwise identifying boundaries of the same length with opposite orientation. Since orientation is a global property, it is sufficient to pick, for instance, an orientation of the first boundary component $\beta_1$.

For a non-orientable surface $K$, Norbury \cite{Nor07} shows a way to make sense of twisting coordinates by choosing a local orientation  in the neighbourhood of $\beta_1$ --- see \cite{PP13} for a different approach. Namely, given a pants decomposition $\mathcal{P}$ on $\Sigma$, cut along curves such that the complement is connected and orientable --- not necessarily along all curves in $\mathcal{P}$. The orientation of $\beta_1$ induces an orientation on such a complement, and this orientation determines twisting coordinates unambiguously. When reconstructing the original non-orientable surface $K$, one identifies two boundaries of the pairs of pants, if they corresponded to a 2-sided curve $\gamma_j$ in $K$, and hence are of the same length $\ell_j$, with a possible twisting $\theta_j$. When such a gluing creates a handle, the identification can be made through a \emph{regular gluing}, i.e.~identifying two boundaries with the opposite orientation, which gives rise to an \emph{orientable handle}, or through an \emph{antipodal gluing}, i.e.~identifying two boundaries with the same orientation, which yields a \emph{non-orientable handle}.
Finally, note that 1-sided curves $\alpha_i$ open M\"obius strips along the core and hence are just glued back to themselves via the antipodal map. As a consequence, no twisting coordinate is attached to 1-sided curves.



For convenience, we henceforth consider a Klein surface $K$ equipped with an orientation on its first boundary component $\beta_1$. When defining the moduli space, care must be taken to ensure that this choice is compatible with the action of the mapping class group. More concretely, let ${\rm MCG}(K)$ be the mapping class group of $K$, which is the group of isotopy classes of self-homeomorphisms fixing each boundary component set-wise. Elements of ${\rm MCG}(K)$ can either preserve or reverse the orientation of $\beta_1$. As shown in \cite{Nor07}, there exists a subgroup  ${\rm MCG}_1(K)\subset {\rm MCG}(K)$ consisting of those mapping classes that preserve the orientation of $\beta_1$.

\begin{lemma}[\cite{Nor07}]\label{lem:orientation}
    For any hyperbolic surface $K$ with at least one boundary, there exists an involution operator $\rho_1\in{\rm MCG}(K)$ which flips the chosen orientation of $\beta_1$. Induced by $\langle\rho_1\rangle$, ${\rm MCG}(K)$ admits an index-two subgroup, denoted by ${\rm MCG}_1(K)$, that preserves the chosen orientation of $\beta_1$.
\end{lemma}

For completeness and clarity, we will include a brief proof of this lemma in the case of stable surfaces, which are the only ones relevant for our purposes, following an analysis of the base cases with Euler characteristic $-1$ (see Lemma~\ref{lem:flipping-or-element}).


\subsubsection{Norbury form}
When $K$ is orientable and oriented, it is well-known that $\mathcal{T}(K)$ admits a symplectic form $\omega_{{\rm WP}}$, called the \emph{Weil--Petersson form}, for which Fenchel--Nielsen coordinates are Darboux coordinates
\begin{equation*}
    \omega_{{\rm WP}}=\sum_{j=1}^{3g-3+n} d\ell_{\gamma_j}\wedge d\theta_j.
\end{equation*}
Although it is written in terms of specific coordinates, it can be shown that the Weil--Petersson form $\omega_{{\rm WP}}$ is invariant under the action of the mapping class group (and independent of the choice of pants decomposition). Therefore, one can define the symplectic volume on the moduli space, as we will discuss shortly.

When $K$ is non-orientable, however, it is evident that there is no symplectic structure because $\dim_\mathbb{R}\mathcal{T}(K)$ can be odd. Nonetheless, Norbury showed in \cite{Nor07} that there exists a top-form on $\mathcal{T}(K)$ with certain properties, which we will call the \emph{Norbury form}:

\begin{definition}\label{def:nu}
    Let $(K,\rho)$ be a non-orientable bordered Klein surface and $\mathcal{T}(K)$ be its Teichm\"{u}ller space. The \emph{Norbury form} $\nu_{{\rm N}}$ in Fenchel--Nielsen coordinates is defined by
    \begin{equation}\label{eq:NorbForm}
\nu_{{\rm N}}\coloneqq\left|\left(\bigwedge_{i=1}^{N_1}\frac{d\ell_{\alpha_i}}{\tanh\frac{\ell_{\alpha_i}}{2}}\right)\wedge\left(\bigwedge_{j=1}^{N_2} d\ell_{\gamma_j}\wedge d\theta_j\right)\right|,
\end{equation}
where the notation $\left|\cdot\right|$ stands for the absolute value of the form.
\end{definition}

\begin{theorem}[\cite{Nor07}]\label{thm:Nor07}
The Norbury form $\nu_{{\rm N}}$ is independent of the choice of pants decompositions, and therefore, it is invariant under ${\rm MCG}(K)$, the mapping class group of $K$.
\end{theorem}

The Norbury form is originally introduced without the absolute value, and is shown to be invariant under the action of the mapping class group up to sign. In Section~\ref{sec:euler characteristic -1 from geometry}, we will explicitly examine actions that would introduce such sign ambiguities. For our purposes, however, it suffices to consider the absolute value of the form, which becomes fully invariant under the mapping class group.


Motivated by Lemma \ref{lem:orientation}, Theorem \ref{thm:Nor07} and the previous discussions, one may define the \emph{moduli space} $\mathcal{M}(K)\coloneqq\mathcal{T}(K)/{\rm MCG}_1(K)$, and consider its volume by taking the Norbury form $\nu_{{\rm N}}$ as the measure\footnote{It is further claimed in \cite{Nor07} that the Norbury form can be defined without taking the absolute value on $\mathcal{T}(K)/{\rm MCG}_1(K)$. This statement needs to be improved because there exists an element in ${\rm MCG}_1(K)$ which preserves the orientation of $\beta_1$ yet flips the sign of the Norbury form (if defined without absolute value) as we will explicitly discuss below. We note nonetheless that Lemma~\ref{lem:orientation} holds.}. When $K$ is non-orientable, however, one encounters an immediate issue: since the Norbury form has singularities whenever $\ell_{\alpha_i}=0$ for any $i$, the integration would be divergent.
We will next present how to resolve this issue following \cite{Gen17}.

\subsubsection{Gendulphe's regularisation}\label{sec:reg}

Gendulphe resolved the singularities of the Norbury form in \cite{Gen17}. The \emph{systole} of $K$ is the length of the shortest simple geodesics of $K$. Let us consider a subset $\mathcal{T}^{\epsilon}(K)\subset \mathcal{T}(K)$ such that the systole of 1-sided geodesics is equal to or greater than $\epsilon\in\mathbb{R}_{>0}$. We choose $\epsilon$ to be small enough in the sense that any two closed geodesics of length at most $\epsilon$ cannot intersect in any hyperbolic surface (independent of the topology).

\begin{remark}\label{rem:small epsilon}
Note that such $\epsilon$ exists and can actually be chosen such that
\begin{equation*}
\sinh\frac{\epsilon}{2}\leq 1.
\end{equation*}
This can be justified directly by the classical collar lemma on (oriented) hyperbolic surfaces (see \cite[Chapter 4]{Bus10}). More precisely, for a non-orientable surface $K$, there exists a set of geodesics~$\Gamma$ such that $K\setminus \Gamma$ is orientable and connected as in Section~\ref{sec:twist coordinates}. Applying the collar lemma to $K\setminus \Gamma$ yields the result also for $K$ non-orientable.

\end{remark}

One may ask how the mapping class group ${\rm MCG}_1(K)$ acts on such a regularised space $\mathcal{T}^{\epsilon}(K)$. The following summarises the results proven by Gendulphe that are most relevant to us:

\begin{theorem}[\cite{Gen17}]\label{thm:Gen17}
For $\epsilon\in\mathbb{R}_{>0}$ small enough as in Remark~\ref{rem:small epsilon}, $\mathcal{T}^{\epsilon}(K)$ is invariant under the mapping class group for any $K$, hence $\mathcal{M}^\epsilon(K)\coloneqq\mathcal{T}^\epsilon(K)/{\rm MCG}_1(K)$ is well-defined. Furthermore, the integral $\int_{\mathcal{M}^\epsilon(K)}\nu_{{\rm N}}$ is finite.
\end{theorem}
From now on we will often employ the following notation, which fixes the topology type and the boundary lengths of our surface $K$.
For $n\in\mathbb{Z}_{\geq1}$ and $g\in\frac12\mathbb{Z}_{\geq0}$ with $2g-2+n>0$, we denote $\mathcal{M}_{g,n}^{\pm,\epsilon}(L_{[n]})$ the \emph{regularised moduli space} of orientable (resp.~non-orientable) hyperbolic surfaces of genus $g$ with $n$ geodesic boundaries $\beta_{[n]}$ of lengths $L_{[n]}$, together with an orientation of $\beta_1$. We denote $\mathcal{T}_{g,n}^{\pm,\epsilon}(L_{[n]})$ the corresponding \emph{regularised Teichm\"uller space}, with $\mathcal{M}_{g,n}^{\pm,\epsilon}(L_{[n]})=\mathcal{T}_{g,n}^{\pm,\epsilon}(L_{[n]})/{\rm MCG}^{\pm}_{g,n}$, where ${\rm MCG}^{\pm}_{g,n}$ refers to the mapping class group of an orientable ($+$) or non-orientable ($-$) surface of topology $(g,n)$, i.e.~the group of isotopy classes of self-homeomorphisms fixing each boundary component set-wise and preserving the orientation of the first boundary $\beta_1$. Note that there is no regularisation effect when $K$ is orientable, and in particular, $\mathcal{M}_{g,n}^{+,\epsilon}(L_{[n]})=\mathcal{M}_{g,n}^{+}(L_{[n]})$, which precisely corresponds to the standard moduli spaces of \emph{oriented} hyperbolic surfaces that have been intensively studied in the literature. Thanks to Theorem~\ref{thm:Nor07} and Theorem~\ref{thm:Gen17}, it is possible to define the volume as follows:
\begin{definition}\label{def:V}
Let $\mathcal{M}_{g,n}^{\pm,\epsilon}(L_{[n]})$ be the regularised moduli space. Its \emph{volume} $V_{g,n}^{\pm,\epsilon}(L_{[n]})$, which we often call the \emph{Gendulphe--Norbury (regularised) volume}, is defined by
\begin{equation*}
V_{g,n}^{\pm,\epsilon}(L_{[n]})\coloneqq\begin{dcases}
 \frac{1}{(3g-3+n)!}\int_{\mathcal{M}_{g,n}^{\pm,\epsilon}(L_{[n]})}\nu_{{\rm N}}\,,   & \text{if } g\in\mathbb{Z}_{\geq0}\,,\\
 \frac{1}{(3g'-2+n)!}\int_{\mathcal{M}_{g,n}^{-,\epsilon}(L_{[n]})}\nu_{{\rm N}}\,,   & \text{if } g=g'+\frac12\in\mathbb{Z}_{\geq0}+\frac12\,.
\end{dcases}
\end{equation*} 
\end{definition}
We chose the normalisation as above for consistency with refined topological recursion, which we will discuss shortly. We leave a more geometric justification (e.g.~an interpretation possibly in terms of contact geometry) to future work.

\subsection{McShane--Norbury identities}\label{sec:MN identity}
Our interest is in finding a recursive structure for $V_{g,n}^{\pm,\epsilon}(L_{[n]})$. In the orientable setting, Mirzakhani utilised McShane--Mirzakhani identities to derive such a recursion \cite{Mir06}. In the non-orientable setting, Norbury proved analogous identities which we will review here. As a preparation, let us introduce four functions which appear in the identities:
\begin{align}
R(x,y,z)=&x-\log\frac{\cosh\frac{y}{2}+\cosh\frac{x+z}{2}}{\cosh\frac{y}{2}+\cosh\frac{x-z}{2}},\label{R(x,y,z)}\\
D(x,y,z)=&R(x,y,z)+R(x,z,y)-x,\label{D(x,y,z)}\\
E(x,y,z)=&R(x,2z,y)-\frac{x}{2},\label{E(x,y,z)}\\
F(x,y,z)=&x-2a(x,y,z)-2a(x,z,y),\label{F(x,y,z)}
\end{align}
where\footnote{The expressions for $a(x,y,z)$ and $a(x,z,y)$ here coincide with $a$ and $b$, respectively, in \cite{Nor07}. Also, there is a misprint in the final expression of $F(x,y,z)$ in \cite{Nor07}, though the derivation is correct.}
\[\tanh a(x,y,z)=\frac{\sinh\frac{x}{2}\sinh z\tanh\frac{z}{2}}{\cosh y+\cosh\frac{x}{2}\cosh z}\,.\]
We note that $R(x,y,z)$ and $D(x,y,z)$ coincide with those in \cite{Mir06}.

When $K$ is a one-bordered Klein bottle, there exists a unique 2-sided primitive geodesic $\gamma$ and there are infinitely many 1-sided geodesics $\alpha_{\mathbb{Z}}$. Given a 1-sided geodesic $\alpha_i$, for some $i\in\mathbb{Z}$, one can always find exactly two 1-sided geodesics not intersecting $\alpha_i$. We label these by $\alpha_{i-1}$ and $\alpha_{i+1}$; note that $\alpha_{i-1}$ and $\alpha_{i+1}$ intersect exactly once. Accordingly, the indexing of the family $\alpha_{\mathbb{Z}}$ is chosen so that $\{\alpha_i,\alpha_{i+1}\}$ is a pair of non-intersecting 1-sided geodesics for every $i\in\mathbb{Z}$.


\begin{theorem}[\cite{Nor07}]\label{prop:Nor07-1}
Let $\mathcal{T}_{1,1}^{-,\epsilon}(L_1)$ be the Teichm\"{u}ller space of one-bordered hyperbolic Klein bottles of geodesic boundary length $L_1$. We denote by $(\ell,\theta)$ the length of the unique 2-sided geodesic $\gamma$ and its twisting coordinate, and by $\ell_i$ the length of the 1-sided geodesic $\alpha_i$ for every $i\in\mathbb{Z}$. Then, on $\mathcal{T}_{1,1}^{-,\epsilon}(L_1)$ the following identity holds:
\begin{equation}
L_1=D(L_1,\ell,\ell)+\sum_{i\in\mathbb{Z}}F(L_1,\ell_i,\ell_{i+1}).\label{MN identity for KB}
\end{equation}
\end{theorem}

When $K$ is a non-orientable surface with $\chi(K)<-1$, it turns out that an analogous formula looks slightly different. Let us denote by $\beta_{[n]}, L_{[n]}$ the set of geodesic boundaries of $K$ and their lengths. We consider a pants decomposition $\mathcal{P}\subset K$ and focus on the pair of pants $P_1$ which has $\beta_1$ as one of its boundaries. Based on the types of the other two boundaries, there are three possible scenarios for $P_1$:
\begin{enumerate}
\item one is a 2-sided curve $\gamma$ and the other one is one of the $\beta_i$, where $i\neq1$,
\item both of them are 2-sided curves $\gamma',\gamma''$,
\item one is a 2-sided curve $\bar\gamma$ and the other one is a 1-sided curve $\alpha$.
\end{enumerate}
Then, the following identity holds:

\begin{theorem}[\cite{Nor07}]\label{prop:Nor07-2}
Let $K$ be an $n$-bordered non-orientable surface of boundary lengths $L_{[n]}$ with $\chi(K)<-1$. On $\mathcal{T}(K)$, we have
\begin{equation}
L_1=\sum_{j=2}^n\sum_{\gamma}R(L_1,L_j,\ell_\gamma)+\sum_{\gamma',\gamma''}D(L_1,\ell_{\gamma'},\ell_{\gamma''})+\sum_{\bar\gamma,\alpha}E(L_1,\ell_{\bar\gamma},\ell_\alpha),\label{MN identity}
\end{equation}
where $\ell_\gamma$ is the geodesic length of a 2-sided curve $\gamma$, and similarly for $\gamma',\gamma'',\bar\gamma, \alpha$, and the summation is taken over all homotopically inequivalent geodesics.
\end{theorem}

We call \eqref{MN identity for KB} and \eqref{MN identity} \emph{McShane--Norbury identities}. We note that the McShane--Norbury identity for $\chi(K)<-1$ resembles the orientable analogue in \cite{Mir06}. In fact, if one removes terms involving 1-sided geodesics, the form of \eqref{MN identity} coincides with that in \cite{Mir06} for orientable surfaces. On the other hand, comparisons are not evident when $\chi=-1$ where the McShane identity for the torus is
\begin{equation*}
L_1=\sum_{\gamma}D(L_1,\ell_{\gamma},\ell_{\gamma}),
\end{equation*}
where the summation is taken over all simple curves not intersecting with each other, which does not agree with \eqref{MN identity for KB} after removing  terms involving 1-sided geodesics.

McShane--Norbury identities are derived in \cite{Nor07} by considering probabilities measuring the possible behaviours of geodesics perpendicular to $\beta_1$ on a pair of pants $P_1$. The four functions $R,D,E,F$ are related to a probability of certain fates. However, it is worth emphasising that $E$ can become negative and one cannot think of $E$ by itself as a probability. The correct interpretation is that in the last summation in \eqref{MN identity}, both $E(L_1,\ell_{\bar\gamma},\ell_\alpha)$ and $E(L_1,\ell_{\bar\gamma},\ell_{\alpha'})$ appear, where $(\alpha,\alpha')$  is a unique pair of intersecting 1-sided curves in the corresponding two-bordered real projective plane, and their sum remains non-negative as explained in \cite{Nor07}.


\subsection{Mirzakhani's integration technique}\label{sec:unfolding}

In the orientable setting, Mirzakhani \cite{Mir06} developed a technique to integrate functions on the moduli space, which she then applied to her generalised versions of McShane's identity to obtain her recursion of the volumes of moduli spaces. We briefly describe the main steps of her result \cite[Theorem 8.1]{Mir06}, which we will then apply to our non-orientable setting in Section \ref{sec:one_bordered_klein_bottles}. To lighten notation, only in this subsection we denote by $\mathcal{T}_{g,n}(L_{[n]})$, $\mathcal{M}_{g,n}(L_{[n]})$, $\mathrm{MCG}_{g,n}$ and $V_{g,n}$, the Teichm\"{u}ller space, moduli space, mapping class group and volume in the classical oriented setting, i.e.~we omit the superscript $+$ from the notation used in the rest of the paper.

Let $\gamma=\sum_{i=1}^kc_i\gamma_i$ be a multi-curve on a closed surface $S_{g,n}$ of genus $g$ with $n$ boundary components $\beta_{[n]}=(\beta_1,\ldots,\beta_n)$ and $X\in \mathcal{T}_{g,n}(L_{[n]})$. For any simple closed curve $\gamma_i$ on $S_{g,n}$, let $[\gamma_i]$ denote the homotopy class of $\gamma_i$ and $\ell_{\gamma_i}(X)$ the hyperbolic length of the geodesic representative of $\gamma_i$ on $X$. The length of a multi-curve $\gamma$ is defined by $\ell_{\gamma}(X)=\sum_{i=1}^kc_i\ell_{\gamma_i}(X)$. For $X\in\mathcal{M}_{g,n}(L_{[n]})$, let
$$
\mathcal{O}_{\gamma}\coloneqq \{[\alpha]\mid \alpha\in \mathrm{MCG}_{g,n} \cdot \gamma\}=\{[\phi(\gamma)]\mid\phi\colon S_{g,n}\rightarrow X \, (\text{marking of } X)\}
$$
be the set of homotopy classes of multi-curves in the $\mathrm{MCG}_{g,n}$-orbit of $\gamma$ on $X$. For any function $f\colon \mathbb{R}_+\rightarrow\mathbb{R}_+$, consider\footnote{For simplicity, we consider $f$ as a function of one variable to state Theorem
~\ref{Mirz-unfolding} in the same form as Mirzakhani's original main unfolding result. However, this discussion extends to any function $F\colon\mathbb{R}^k\to\mathbb{R}$, and the unfolding arguments work in the same way (see \cite[Lemma 8.4]{Mir06}). Note that \cite[Lemma 8.4]{Mir06} is stated for multi-curves $\Gamma=(\Gamma_1,\ldots,\Gamma_k)$ given as ordered tuples; if instead one works with unordered multi-curves, one must account for the corresponding symmetry factors, as discussed in Step 1 below.}
$$
\begin{array}{lrcl}
f_{\gamma}\colon & \mathcal{M}_{g,n}(L_{[n]}) &\longrightarrow  & \mathbb{R}_+\\
& X & \mapsto & \sum_{[\alpha]\in\mathcal{O}_{\gamma}}f(\ell_{\alpha}(X)).
\end{array}
$$

Our goal is to be able to integrate functions $f_{\gamma}$ over $\mathcal{M}_{g,n}(L_{[n]})$ with respect to the Weil--Petersson volume form, which is a priori a very difficult task. Mirzakhani's idea consists in unfolding the integral to a cover  $\mathcal{M}_{g,n}^{\gamma}(L_{[n]})$ of the moduli space, over which it will be simpler to integrate, since  $\mathcal{M}_{g,n}^{\gamma}(L_{[n]})$ can be described in terms of moduli spaces of simpler surfaces, obtained after cutting along $\gamma$. Consider the covering $\pi^{\gamma}\colon \mathcal{M}_{g,n}^{\gamma}(L_{[n]})\rightarrow \mathcal{M}_{g,n}(L_{[n]})$, where
$$
\mathcal{M}_{g,n}^{\gamma}(L_{[n]})\coloneqq \{(X,\alpha)\mid X\in \mathcal{M}_{g,n}(L_{[n]}), \alpha=(\alpha_1,\ldots,\alpha_k), \alpha_i \text{ geodesic with } [\alpha_i]\in \mathcal{O}_{\gamma_i}\}
$$
and $\pi^{\gamma}(X,\alpha)=X$, and the length function 
$$
\begin{array}{lrcl}
\ell \colon & \mathcal{M}_{g,n}^{\gamma}(L_{[n]})&\longrightarrow  & \mathbb{R}_+\\
& Y=(X,\alpha) & \mapsto & \ell(Y)=\ell_{\alpha}(X).
\end{array}
$$
Observe that $\mathcal{M}_{g,n}^{\gamma}(L_{[n]})=\mathcal{T}_{g,n}(L_{[n]})/G_{\gamma}$, where $G_{\gamma}\coloneqq\bigcap_{i=1}^k\mathrm{Stab}(\gamma_i)\subset \mathrm{Stab}(\gamma)\subset\mathrm{MCG}_{g,n}$.

\begin{description}
\item[{\rm Step 1}] \phantomsection \label{item:Step1} For a function $f\circ \ell\colon \mathcal{M}_{g,n}^{\gamma}(L_{[n]})\rightarrow \mathbb{R}_+$, one can consider the push-forward function 
$$
\begin{array}{lrcl}
\pi_*^{\gamma}(f\circ\ell) \colon &  \mathcal{M}_{g,n}(L_{[n]})&\longrightarrow  & \mathbb{R}_+\\
& X & \mapsto & \sum_{Y\in(\pi^{\gamma})^{-1}(X)} f(\ell(Y)).
\end{array}
$$
The first step is to realise that $\pi_*^{\gamma}(f\circ\ell)$ is very closely related to $f_{\gamma}$:
\begin{align*}
\pi_*^{\gamma}(f\circ\ell)(X) & =\sum_{h\in \mathrm{MCG}_{g,n}/G_{\gamma}} f(\ell_{h\cdot\gamma}(X)) \nonumber \\
&=\left|\mathrm{Sym}(\gamma)\right|\sum_{h\in\mathrm{MCG}_{g,n}/\mathrm{Stab}(\gamma)}f(\ell_{h\cdot\gamma}(X))=\left|\mathrm{Sym}(\gamma)\right| f_{\gamma}(X),
\end{align*}
where $\mathrm{Sym}(\gamma)\coloneqq \mathrm{Stab}(\gamma)/G_{\gamma}$ is the symmetry group of $\gamma$. Observe that $\left|\mathrm{Sym}(\gamma)\right|\neq 1$ will impose that $c_i=c_j$, if $h(\gamma_i)=\gamma_j$ for $h\in\mathrm{Sym}(\gamma)$; therefore, for $h\in\mathrm{Sym}(\gamma)$
$$
\ell_{h\cdot\gamma}(X)=\sum_{i=1}^kc_i\ell_{h\cdot\gamma_i}(X)=\sum_{j=1}^kc_j\ell_{\gamma_j}(X)=\ell_{\gamma}(X).
$$

\item[{\rm Step 2}] \phantomsection \label{item:Step2} The second step, which best captures Mirzakhani's idea, is to observe that the fibres of $\pi^{\gamma}$ correspond precisely to the sets over which we sum in the definition of $\pi_*^{\gamma}(f\circ\ell)$. Thus
\begin{equation*}
\int_{\mathcal{M}_{g,n}(L_{[n]})}\pi_*^{\gamma}(f\circ\ell)(X)\,dX=\int_{\mathcal{M}_{g,n}^{\gamma}(L_{[n]})} f\circ\ell(Y)\,dY,
\end{equation*}
where $dX=\frac{\omega_{{\rm WP}}^{3g-3+n}}{(3g-3+n)!}$ and $dY=(\pi^{\gamma})^*dX$.


\item[{\rm Step 3}] \phantomsection \label{item:Step3} The third step utilises the fact that the covering space $\mathcal{M}_{g,n}^{\gamma}(L_{[n]})$ is closely related to the moduli space of hyperbolic structures on $S_{g,n}(\gamma)$. Let $\mathcal{M}(S_{g,n},\ell_{\beta}=L_{[n]})\coloneqq \mathcal{M}_{g,n}(L_{[n]})$ and $S_{g,n}(\gamma)$ be the (possibly disconnected) surface cut along $\gamma$, with $n+2k$ boundaries and $s=s(\gamma)$ connected components $S_i$ of genus $g_i$, each with an ordered tuple of boundary components denoted $\partial S_i$ of lengths $\ell_{\partial S_i}$. Then, we have: 
\begin{equation*}
\int_{\mathcal{M}_{g,n}^{\gamma}(L_{[n]})} f(\ell(Y))\,dY
=\int_{\underline{l}\in\mathbb{R}_+^k}\int_{\theta_i\in[0,\ell_i]}f(|\underline{l}|)\,V(\mathcal{M}(S_{g,n}(\gamma),\ell_{\beta}=L_{[n]},\ell_{\gamma}=\underline{l}))\,dl_1d\theta_1\cdots dl_kd\theta_k,
\end{equation*}
where $\underline{l}=(l_1,\ldots,l_k)$, $|\underline{l}|=c_1l_1+\cdots + c_kl_k$, $V(\mathcal{M}(S_i,\ell_{\partial S_i}))\coloneqq V_{g_i,|\partial S_i|}(\ell_{\partial S_i})$ and
\[V(\mathcal{M}(S_{g,n}(\gamma),\ell_{\beta}=L_{[n]},\ell_{\gamma}=\underline{l}))\coloneqq\prod_{i=1}^sV(\mathcal{M}(S_i,\ell_{\partial S_i})).\]

\end{description}

Putting together these three outlined steps, one obtains Mirzakhani's result:
\begin{theorem}[{\cite[Theorem 8.1]{Mir06}}]\label{Mirz-unfolding}
For any multi-curve $\gamma=\sum_{i=1}^kc_i\gamma_i$, the integral of $f_{\gamma}$ over $\mathcal{M}_{g,n}(L_{[n]})$, for $2g-2+n>0$, with respect to the Weil--Petersson volume form is given by\footnote{Observe that in Mirzakhani's result the pre-factor appearing includes a $2$ for every component $\gamma_i$ of the multi-curve separating off a one-handle. This discrepancy is due to the fact that Mirzakhani ignores the elliptic involution on $\mathcal{M}_{1,1}(L_{1})$.}
\begin{equation*}
\int_{\mathcal{M}_{g,n}(L_{[n]})} f_{\gamma}(X)\, dX
=\frac{1}{|\mathrm{Sym}(\gamma)|}\int_{\underline{l}\in\mathbb{R}_+^k}l_1\cdots l_k\, f(|\underline{l}|)\, V(\mathcal{M}(S_{g,n}(\gamma),\ell_{\beta}=L_{[n]},\ell_{\gamma}=\underline{l}))\, dl_1\cdots dl_k.
\end{equation*}
\end{theorem}

To calculate the volumes $V_{g,n}(L_{[n]})$ using this technique, Mirzakhani expressed the constant function on $\mathcal{M}_{g,n}(L_{[n]})$ as the sum of functions of the form of $f_{\gamma}$, generalising McShane's original identity on a punctured torus to any surface with boundaries.

\section{Euler characteristic $-1$ from hyperbolic geometry}\label{sec:euler characteristic -1 from geometry}

For all $(g,n)$ with $2-2g-n\geq0$, we set $V^{\pm,\epsilon}_{g,n}(L_{[n]})=0$ because the corresponding moduli space is not well-defined\footnote{Some literature (e.g.~\cite{EO07-2,Sta23}) conventionally sets unstable volumes to some functions, but they do not have any geometric meaning.}. There are four topologically different bordered Klein surfaces with $\chi(K)=-1$: a three-bordered sphere, a one-bordered torus, a two-bordered real projective plane, and a one-bordered Klein bottle. The volumes for the first two cases are well-known to be
\begin{equation}
V_{0,3}^+(L_1,L_2,L_3)=1,\quad V_{1,1}^+(L_1)=\frac{L_1^2}{48}+\frac{\pi^2}{12}.\label{V_{0,3} and V_{1,1}}
\end{equation}
Note that $V_{1,1}^+(L_1)$ is indeed half of the volume computed in \cite{Mir06}. This is because we have already taken care of the extra symmetry existing only for $\mathcal{M}_{1,1}^+(L_1)$, i.e. the elliptic involution.

\subsection{Two-bordered real projective planes}\label{sec:1/2,2}
When $K$ is a two-bordered real projective plane, let us fix the boundary lengths to be  $(L_1,L_2)$. Up to homotopy, there are exactly two 1-sided geodesics $\alpha,\alpha'$ which intersect once, and we denote their lengths by $\ell_{\alpha},\ell_{\alpha'}$, respectively. In the hyperbolic setting, these lengths are mutually related by \cite{Nor07}
\begin{equation}
\cosh\frac{L_1}{2}+\cosh\frac{L_2}{2}=2\sinh\frac{\ell_\alpha}{2}\sinh\frac{\ell_{\alpha'}}{2}.\label{constraint for Mobius}
\end{equation}
That is, if we fix one of $\ell_{\alpha},\ell_{\alpha'}$, then the other one is uniquely fixed. This is consistent with the fact that $\mathcal{T}(K)\cong\mathbb{R}_{>0}$ on which either of $\ell_{\alpha},\ell_{\alpha'}$ can be taken as a coordinate. It can be easily deduced from \eqref{constraint for Mobius} that the sign of the Norbury form is flipped under the exchange of $\alpha$ and $\alpha'$.

\begin{figure}[ht]
    \centering
    \includegraphics[width=0.95\textwidth]{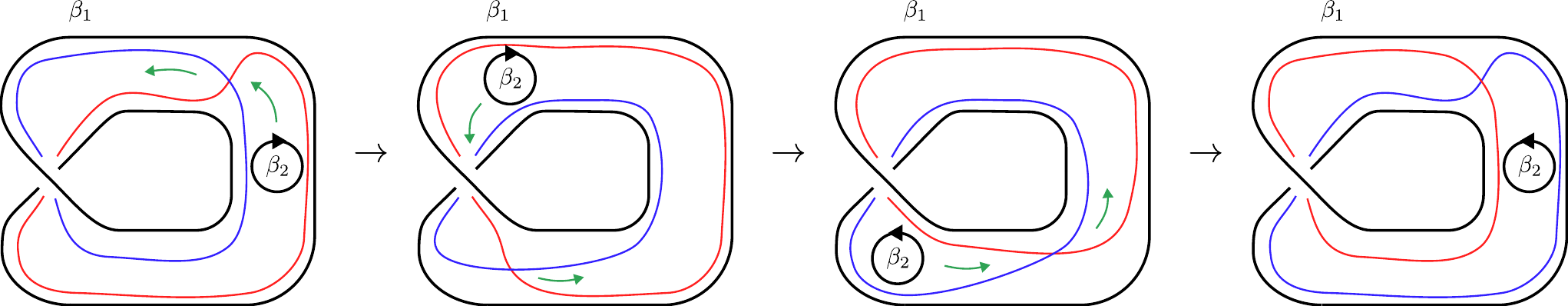}
    \caption{The picture shows the effect of an isotopy $h_t$ of $K$ which moves the second boundary $\beta_2$ around the M\"obius strip once, while keeping the boundary $\beta_1$ (of the M\"{o}bius strip) fixed. We have $h_0=\mathrm{Id}$, $\left.h_t\right|_{\beta_1}=\mathrm{Id}$ and $h_1=Y^{(1)}$, which is still the identity on $\beta_1$ and reverses the orientation on $\beta_2$.}\label{fig:Y1} 
\end{figure}

The mapping class group of $K$ has two generators $Y^{(1)}, Y^{(2)}$ that exchange $\alpha$ and $\alpha'$ \cite{Kor02,Gen17}. The homeomorphism $Y^{(i)}$ is the identity on the $i$th boundary and reverses the orientation of the other boundary (see Figure~\ref{fig:Y1}). These involutions are instances of the so-called \emph{$Y$-homeomorphisms}; in general, mapping class groups of non-orientable surfaces are generated not only by Dehn twists but also by $Y$-homeomorphisms \cite{Lic63,PP13}. Since we are interested in the index-two subgroup ${\rm MCG}^1(K)$ of the mapping class group that preserves the orientation of the first boundary, we do not consider the second generator $Y^{(2)}$:
\begin{align*}
{\rm MCG}^1(K)&={\rm MCG}^{-}_{\frac12,2}=\langle Y^{(1)}\rangle\\
 \mathcal{M}^{-,\epsilon}_{\frac12,2}(L_1,L_2)&=\mathcal{T}^{-,\epsilon}_{\frac12,2}(L_1,L_2)/{\rm MCG}^{-}_{\frac12,2}.
\end{align*}


\begin{proposition}\label{prop:V_{1/2,2}}
For a two-bordered real projective plane of boundary lengths $L_{[2]}$, we have\footnote{This proof is also outlined in \cite{Sta23}.}
\begin{equation}\label{V_{1/2,2}-1}
V^{-,\epsilon}_{\frac12,2}(L_{[2]})=\log\frac{\cosh\frac{L_1+L_2}{4}\cosh\frac{L_1-L_2}{4}}{\sinh^2\frac{\epsilon}{2}}.
\end{equation}
\end{proposition}
\begin{proof}
Since the mapping class group exchanges $\alpha$ with $\alpha'$, we can assume without loss of generality that $\ell_\alpha\leq\ell_{\alpha'}$ in the moduli space. Let $\ell_*$ be defined by the point in which $\ell_\alpha=\ell_{\alpha'}$, i.e. by
\begin{equation*}
    2\cosh\frac{L_1+L_2}{4}\cosh\frac{L_1-L_2}{4}=\cosh\frac{L_1}{2}+\cosh\frac{L_2}{2}=2\sinh^2\frac{\ell_*}{2}.
\end{equation*}
Then, a fundamental (oriented) domain $\mathcal{F}_{\frac12,2}(L_{[2]})$ of the moduli space can be realised as $\mathcal{F}_{\frac12,2}(L_{[2]})=\{\ell_\alpha\in\mathbb{R}_{>0}\mid 0<\ell_\alpha\leq\ell_*\}$. Accordingly, we consider the regularised fundamental domain $\mathcal{F}_{\frac12,2}^\epsilon(L_{[2]})=\{\ell_\alpha\in\mathbb{R}_{>0}\mid\epsilon<\ell_\alpha\leq\ell_*\}$ endowed with the measure $\nu_{{\rm N}}=\frac{d\ell_\alpha}{\tanh\frac{\ell_\alpha}{2}}$ which remains positive and finite within the regularised fundamental domain. Once we have the correct integral domain and measure, the regularised volume can be explicitly computed in terms of elementary functions: 
\begin{align}
V^{-,\epsilon}_{\frac12,2}(L_{[2]}) & \coloneqq\int_{\mathcal{M}_{\frac12,2}^{-,\epsilon}(L_{[2]})}\left|\frac{d\ell_\alpha}{\tanh\frac{\ell_\alpha}{2}}\right|=\int_\epsilon^{\ell_*}\frac{d\ell_\alpha}{\tanh\frac{\ell_\alpha}{2}}\nonumber \\
&=2\log\sinh\frac{\ell_*}{2}-2\log\sinh\frac{\epsilon}{2}.\nonumber
\end{align}
We obtain the proposition by substituting the defining equation of $\ell_*$.
\end{proof}

Although this suffices to prove Proposition~\ref{prop:V_{1/2,2}}, let us yet give a slightly different derivation. Since the mapping class group is of order 2, one can integrate over the Teichm\"{u}ller space with a factor of one-half in front:
\begin{align}
V^{-,\epsilon}_{\frac12,2}(L_{[2]})=\frac12\int_{\mathcal{T}_{\frac12,2}^{-,\epsilon}(L_{[2]})}\frac{d\ell_\alpha}{\tanh\frac{\ell_\alpha}{2}}.\label{V_{1/2,2}-2}
\end{align}
Perhaps counterintuitively, we note that $\mathcal{T}_{\frac12,2}^{-,\epsilon}(L_{[2]})\not\cong\mathbb{R}_{>\epsilon}$. This is because the constraint \eqref{constraint for Mobius} gives the $\epsilon$-dependent \emph{upper bound} $\Lambda$ determined by
\begin{equation}
\cosh\frac{L_1}{2}+\cosh\frac{L_2}{2}=2\sinh\frac{\epsilon}{2}\sinh\frac{\Lambda(L_1,L_2,\epsilon)}{2}.\label{Lambda}
\end{equation}
Thus, the correct regularised Teichm\"{u}ller space is 
\begin{equation}\label{Teich12-2}
\mathcal{T}_{\frac12,2}^{-,\epsilon}(L_{[2]})\cong\{\ell_\alpha\in\mathbb{R}_{>0}\mid\epsilon<\ell_\alpha<\Lambda(L_1,L_2,\epsilon)\}.
\end{equation}
One can easily check that the integral \eqref{V_{1/2,2}-2} agrees with \eqref{V_{1/2,2}-1}. It is crucial to keep in mind that Gendulphe's regularisation imposes not only the lower bound for 1-sided geodesics but also the upper bound, as evidently seen in this example.

\subsection{One-bordered Klein bottles}\label{sec:one_bordered_klein_bottles}





When $K$ is a one-bordered Klein bottle, recall that it contains a unique 2-sided primitive geodesic~$\gamma$ and infinitely many 1-sided geodesics $(\alpha_i)_{i\in\mathbb{Z}}$, indexed so that $\alpha_i\cap\alpha_{i+1}=\varnothing$ for every $i\in\mathbb{Z}$. The mapping class group $\mathrm{MCG}(K)$ is isomorphic to $\mathbb{D}_{\infty}\times\mathbb{Z}/2\mathbb{Z}$, where $\mathbb{D}_{\infty}$ is the infinite dihedral group \cite{Kor02,Sze06,Gen17}. The index two subgroup $\mathbb{Z}\subset \mathbb{D}_{\infty}$ is generated by the Dehn twist $D_\gamma$ along $\gamma$, whose action is $D_\gamma:\alpha_i\mapsto\alpha_{i+1}$ for all $i\in\mathbb{Z}$. In order to describe the other generators, which are certain $Y$-homeomorphisms, let us view a one-bordered Klein bottle with boundary $\beta_1$ as a connected sum of a two-bordered projective plane $N_{\frac12,2}$ and a M\"{o}bius strip $M$ such that the first boundary of $N_{\frac12,2}$ is $\beta_1$ whereas the other boundary $\beta_{\mathrm{int}}$ is glued to the M\"{o}bius strip $M$. 

For $i=1,2$, let $\widetilde{Y}^{(i)}$ be the $Y$-homeomorphism such that 
$$
\left.\widetilde{Y}^{(i)}\right|_{N_{\frac12,2}}=Y^{(i)},
$$ 
where $Y^{(i)}$ are those in Section~\ref{sec:1/2,2}. That is, $\widetilde{Y}^{(1)}$ ($\widetilde{Y}^{(2)}$) fixes (reverses) the orientation of the boundary $\beta_1$, while it reverses (preserves) the orientation of the internal $\beta_{\mathrm{int}}$. Hence $\mathrm{MCG}(K)\cong\langle D_{\gamma}, \widetilde{Y}^{(1)}, \widetilde{Y}^{(2)}\rangle$, while
$$
\mathrm{MCG}^1(K)=\mathrm{MCG}_{1,1}^-\cong\langle D_{\gamma}, \widetilde{Y}^{(1)}\rangle.
$$
The action of $\widetilde{Y}^{(1)}=Y_0$ is to send $\alpha_{j}\mapsto\alpha_{-j}$, after having chosen a reference 1-sided geodesic $\alpha_0$. One can consider $Y_i=D_\gamma^{2i}\circ Y_0$ for each $i\in\mathbb{Z}$, which fixes $\alpha_i$ and sends $\alpha_{i+j}\mapsto\alpha_{i-j}$ for all $j\in\mathbb{Z}$. Note that the unique 2-sided $\gamma$ is invariant under both $D_\gamma$ and $Y_i$ as expected.

%

In terms of coordinates, let us consider $s_i\coloneqq\sinh\frac{\ell_{i}}{2}$ where $\ell_i=\ell_{\alpha_i}$. The Norbury form is then given by\footnote{The coordinate $s$ in the present paper is equivalent to $Y$ in \cite{Nor07}. However, note that there is a typo in Eq.~(25) in \cite{Nor07} --- the overall factor of 4 is missing.}
\begin{equation*}
    \nu_{{\rm N}}=4\frac{ds_i}{s_i}\wedge\frac{ds_{i+1}}{s_{i+1}}.
\end{equation*}
Furthermore, one can show that with respect to the Fenchel--Nielsen coordinates associated to the unique 2-sided geodesics $(\ell,\theta)$, every $s_i$ is written as\footnote{To the best of the authors' knowledge, this is not written anywhere explicitly but it is easy to derive by following the trace identity discussed in \cite[Section 3]{Nor07}.}
\begin{equation*}
    s_i=\frac{\cosh\frac{L_1}{4}\cosh\frac{\theta+i\ell}{2}}{\sinh\frac{\ell}{2}}, \text{ for } i\in\mathbb{Z}.
\end{equation*}
The action of the Dehn twist and the $Y$-homeomorphism are respectively given by $D_\gamma:(\ell,\theta)\mapsto(\ell,\theta+\ell)$ and $Y_0:(\ell,\theta)\mapsto(\ell,-\theta)$. Thus, the Norbury form $\nu_{{\rm N}}$ is invariant under $D_\gamma$ and under each $Y_i$ (while it would be only anti-invariant under each $Y_i$ if it had been defined without an absolute value, as in \cite{Nor07}). 

We are now ready to give a brief proof of the following lemma, which implies Lemma~\ref{lem:orientation} for stable surfaces as a corollary.

\begin{lemma}\label{lem:flipping-or-element}
On any stable surface $S$ with at least one boundary, there exists a homeomorphism that flips a chosen orientation on a boundary component.
\end{lemma}
\begin{proof}
Every pair of pants possesses an orientation-reversing reflection. This can be extended to any orientable surface $S$, since it can be decomposed into pairs of pants.

Now let $S$ be non-orientable. If $\chi(S)=-1$, that is, $S$ is a two-bordered projective plane (resp.~a one-bordered Klein bottle), then $Y^{(2)}$ (resp.~$\widetilde{Y}^{(2)}$) is the homeomorphism that flips the orientation of the first boundary component.
If $\chi(S)<-1$, then $S$ is the connected sum of a two-bordered projective plane $N$ and a stable surface $\widetilde{S}$. We know that $Y^{(2)}$ flips the orientation of the first boundary of $N$, that is, the first boundary of $S$, and fixes the orientation of the second boundary of $N$; hence $Y^{(2)}$ can be extended to a homeomorphism of the whole $S$ by the identity.
\end{proof}



We now turn to computing $V_{1,1}^{-,\epsilon}(L_1)$. The first step is to simplify the definition of the volume into a computable form by utilising the McShane--Norbury identity:
\begin{lemma}
The volume $V_{1,1}^{-,\epsilon}(L_1)$ can be expressed as an integral over the regularised Teichm\"{u}ller space $\mathcal{T}_{1,1}^{-,\epsilon}(L_1)$, for any $i\in\mathbb{Z}$, as
\begin{equation}
V_{1,1}^{-,\epsilon}(L_1)=\frac12\int_{\mathcal{T}_{1,1}^{-,\epsilon}(L_1)}\frac{F(L_1,\ell_i,\ell_{i+1})}{L_1-D(L_1,\ell,\ell)}\nu_{{\rm N}},\label{unfolded volume}
\end{equation}
where $D(x,y,z),F(x,y,z)$ are in \eqref{D(x,y,z)} and \eqref{F(x,y,z)}.
\end{lemma}

\begin{proof}
Note that the McShane--Norbury identity for Klein bottles \eqref{MN identity for KB} can also be written as
\begin{equation}
  1=  \frac{\sum_{i\in\mathbb{Z}}F(L_1,\ell_i,\ell_{i+1})}{L_1-D(L_1,\ell,\ell)}.\label{MN identity for KB2}
\end{equation}
Inserting this to the defintion of the Gendulphe--Norbury volume, one finds:
\begin{equation*}
V_{1,1}^{-,\epsilon}(L_1)=\int_{\mathcal{M}_{1,1}^{-,\epsilon}(L_1)}\frac{\sum_{i\in\mathbb{Z}}F(L_1,\ell_i,\ell_{i+1})}{L_1-D(L_1,\ell,\ell)}\nu_{{\rm N}}.
\end{equation*}
Notice that $D_\gamma$ keeps $\ell$ invariant and sends $\ell_i$ to $\ell_{i+1}$ for all $i$, hence in particular, the integrand is invariant under the action of $D_\gamma$ and of each $Y_i$. Applying an analogous technique to Mirzakhani's unfolding from Section~\ref{sec:unfolding} (only \hyperref[item:Step1]{Step 1} and \hyperref[item:Step2]{Step 2} are necessary in this case), we can get rid of the infinite summation in $i\in\mathbb{Z}$ by lifting the integral domain from $\mathcal{M}_{1,1}^{-,\epsilon}(L_1)$ to $\mathcal{T}_{1,1}^{-,\epsilon}(L_1)/G_{\alpha}$, where $\alpha$ is the 1-sided multi-curve $\alpha_i+\alpha_{i+1}$. Observe that $\mathrm{Stab}(\alpha)=\langle D^{2i+1}\circ Y_0\rangle$, hence $|\mathrm{Sym}(\alpha)|=2$, and we obtain the lemma. 
\end{proof}

A few remarks are in order. First, notice that the integral domain becomes the Teichm\"{u}ller space itself. This stands in clear contrast to the case of the torus in \cite{Mir06}, whose integral domain after unfolding is $\mathcal{T}_{1,1}^+(L_1)/\Gamma$, with $\Gamma$ a certain subgroup of the mapping class group of the torus. Next, the integral domain $\mathcal{T}_{1,1}^+(L_1)/\Gamma$ is simple in the case of the torus --- in the standard coordinates $(\ell,\theta)$ it is $0\leq\ell<\infty$ and $0\leq\theta<\ell$. On the other hand, $\mathcal{T}_{1,1}^{-,\epsilon}(L_1)$ is complicated due to the $\epsilon$-regularisation, and hence we need a more careful analysis to derive $V_{1,1}^{-,\epsilon}(L_1)$ --- imposing the lower bound of 1-sided geodesics may also impose an upper bound as we saw for the computation of $V_{\frac12,2}^{-,\epsilon}$.

With these points noted, the following theorem shows that, by accounting for such nontrivial regularisation effects, one can explicitly compute $V_{1,1}^{-,\epsilon}$, which has a remarkably simple expression.

\begin{theorem}\label{thm:KB}
\begin{equation}
V_{1,1}^{-,\epsilon}(L_1)=-{\rm Li}_2\left(-\frac{\cosh^2\frac{L_1}{4}}{\sinh^2\frac{\epsilon}{2}}\right).
\end{equation}
\end{theorem}
\begin{proof}
As shown in \cite{Nor07} we have for any $i\in\mathbb{Z}$ that
\begin{equation}
s_i^2+s_{i+1}^2-2s_is_{i+1}\cosh\frac{\ell}{2}=-\cosh^2\frac{L_1}{4},\quad s_{i-1}s_{i+1}=s_i^2+\cosh^2\frac{L_1}{4},\label{trace identity}
\end{equation}
It follows from \eqref{trace identity} that if $s_{i-1}\geq s_{i}$ for some $i$, then we have  $\cdots\geq s_{i-2}\geq s_{i-1}\geq s_{i}$, and similarly if $s_{i}\leq s_{i+1}$ for some $i$, then it follows that  $s_{i}\leq s_{i+1}\leq s_{i+2}\leq\cdots$. However, $s_{i}\leq s_{i+1}$ does not impose any inequality between $s_{i-1}$ and $s_{i}$. Therefore, we would like to decompose $\mathcal{T}_{1,1}^{-,\epsilon}(L_1)$ into $\bigcup_i U_i$ such that $s_i$ is smaller than the other $s_{j\neq i,j\in\mathbb{Z}}$ in every $U_i$.

Fix $i\in\mathbb{Z}$ and assume that $s_i\leq s_{i+1}$. We would like to find the condition on $s_i,s_{i+1}$ such that $s_{i-1}\geq s_i$. Such a condition can be determined by the second equation in \eqref{trace identity} by substituting $s_{i-1}=s_i$. More concretely, when $s_i,s_{i+1}$ are in the following domain
\begin{equation*}
U_i=\left\{(s_i,s_{i+1})\in\mathbb{R}^2\,\;\Big|\;\,\sinh\frac{\epsilon}{2}\leq s_i,\;\; s_i\leq s_{i+1}\leq \frac{s_i^2+\cosh^2\frac{L_1}{4}}{s_i}\right\},
\end{equation*}
then $s_i$ is smaller than any other $s_{j\neq i}$ for all $j\in\mathbb{Z}$. By construction, we have $\mathcal{T}_{1,1}^{-,\epsilon}(L_1)=\bigcup_i U_i$.

We can then simplify the integral \eqref{unfolded volume}. Without loss of generality, we can set $i=0$ in \eqref{unfolded volume}. Since the Norbury form $\nu_{{\rm N}}$ as well as the length of the unique 2-sided geodesic $\ell$ are invariant under the Dehn twist, one notices that the integration of any top form $\omega(\ell_0,\ell_1;\ell)$ on $U_k$ is the same as the integration of $\omega(\ell_{-k},\ell_{-k+1};\ell)$ on $U_0$. In particular,
\begin{equation*}
\int_{U_{k}}\frac{F(L_1,\ell_0,\ell_{1})}{L_1-D(L_1,\ell,\ell)}\nu_{{\rm N}}=\int_{U_{0}}\frac{F(L_1,\ell_{-k},\ell_{1-k})}{L_1-D(L_1,\ell,\ell)}\nu_{{\rm N}}.
\end{equation*}
Therefore, since the intersection of the $U_k$s is of measure zero, we find
\begin{align}
V_{1,1}^{-,\epsilon}(L_1)=&\frac12\sum_{k\in\mathbb{Z}}\int_{U_{k}}\frac{F(L_1,\ell_0,\ell_{1})}{L_1-D(L_1,\ell,\ell)}\nu_{{\rm N}}\nonumber\\
=&\frac12\int_{U_{0}}\frac{\sum_{k\in\mathbb{Z}}F(L_1,\ell_{-k},\ell_{1-k})}{L_1-D(L_1,\ell,\ell)}\nu_{{\rm N}}\nonumber\\
=&2\int_{\sinh\frac{\epsilon}{2}}^{\infty}\frac{ds_0}{s_0}\int_{s_0}^{\frac{s_0^2+\cosh^2\frac{L_1}{4}}{s_0}}\frac{ds_1}{s_1},\nonumber
\end{align}
where the second equality follows from the Fubini--Torelli theorem, since $F(L_1,\ell_{-k},\ell_{1-k})>0$, and the third is obtained by applying the McShane--Norbury identity for Klein bottles, namely Equation \eqref{MN identity for KB2}. The remaining integral can be evaluated explicitly, completing the proof.
\end{proof}

We note that for $z\in\mathbb{C}\backslash\mathbb{R}_{\geq1}$, the dilogarithm satisfies the following reflection property:
\begin{equation}
    {\rm Li}_2(z)+ {\rm Li_2}(z^{-1})=-\frac{\pi^2}{6}-\frac12\log^2(-z).\label{dilog identity}
\end{equation}
By substituting $z = -\frac{\sinh^2\frac{\epsilon}{2}}{\cosh^2\frac{L_1}{4}}\in\mathbb{C}\backslash\mathbb{R}_{\geq1}$, we obtain
\begin{equation}
  V_{1,1}^{-,\epsilon}(L_1)=  2\log^2\frac{\sinh\frac{\epsilon}{2}}{\cosh\frac{L_1}{4}}+\frac{\pi^2}{6}+{\rm Li}_2\left(-\frac{\sinh^2\frac{\epsilon}{2}}{\cosh^2\frac{L_1}{4}}\right).\label{expansion of V_{KB}}
\end{equation}
This explicitly shows that $V_{1,1}^{-,\epsilon}(L_1)$ in fact takes positive values for any $L_1\in\mathbb{R}_{>0}$. In addition, it turns out that this expression is more convenient to compare with refined topological recursion. In particular, the expression \eqref{expansion of V_{KB}} admits the following expansion:

\begin{proposition}\label{prop:expansion of V_KB} For $L_1>2\epsilon$\footnote{The condition $L_1>2\epsilon$ is necessary only for expanding as in the proposition, but of course $L_1$ can be set to any value, e.g. $L_1<2\epsilon$ in $V_{1,1}^{-,\epsilon}(L_1)$.}, one can expand $V_{1,1}^{-,\epsilon}(L_1)$ as follows:
\begin{align*}
        V_{1,1}^{-,\epsilon}(L_1)=&\frac{1}{8}L_1^2-L_1\log\left(2\sinh\frac{\epsilon}{2}\right)+\frac{\pi^2}{6}+2\log^2\left(2\sinh\frac{\epsilon}{2}\right)\nonumber\\
        &-L_1\sum_{k\geq1}\frac{(-1)^k}{k}e^{-\frac{k}{2}L_1}+\sum_{k\geq1}(-1)^k U^\epsilon_ke^{-\frac{k}{2}L_1},
    \end{align*}
    where
    \begin{align}
        U^\epsilon_k=\frac{4}{k}\log\left(2\sinh\frac{\epsilon}{2}\right)+\frac{4}{k}\sum_{j=1}^{k-1}\frac{1}{j}+\sum_{j=1}^k\frac{\left(2\sinh\frac{\epsilon}{2}\right)^{2j}}{j^2}\binom{k+j-1}{2j-1}.\label{U for 1,1}
    \end{align}    
\end{proposition}
    
\begin{proof}
    First, recall that ${\rm Li}_2(z)$ admits the following expansion when $|z|<1$:
    \begin{equation*}
    {\rm Li}_2(z) = \sum_{k\geq1}\frac{z^k}{k^2}.
    \end{equation*}
    Next, we have
    \begin{align*}
        \log^2\left(1+e^{-\frac{L}{2}}\right)=\sum_{i,j\geq1}\frac{(-1)^{i+j}}{ij}e^{-\frac{i+j}{2}L}=\sum_{k\geq2}(-1)^ke^{-\frac{k}{2}L}\sum_{j=1}^{k-1}\frac{1}{j(k-j)}=\sum_{k\geq1}(-1)^ke^{-\frac{k}{2}L}\frac{2}{k}\sum_{j=1}^{k-1}\frac{1}{j}.
    \end{align*}
    Also, note that
    \begin{align*}
        \frac{1}{\left(2\cosh\frac{L}{4}\right)^{2j}}=&e^{-\frac{j}{2}L}\sum_{l\geq0}(-1)^l\binom{2j+l-1}{l}e^{-\frac{l}{2}L},
    \end{align*}
    which follows by differentiating the geometric series $2j-1$ times. It then follows that
    \begin{align*}
        \sum_{j\geq1}\frac{(-1)^j}{j^2}\left(\frac{2\sinh\frac{\epsilon}{2}}{2\cosh\frac{L}{4}}\right)^{2j}=\sum_{k\geq1}(-1)^ke^{-\frac{k}{2}L}\sum_{j=1}^{k}\frac{\left(2\sinh\frac{\epsilon}{2}\right)^{2j}}{j^2}\binom{k+j-1}{2j-1}.
    \end{align*}
The rearrangement of the summation at the equality holds when $L_1>2\epsilon$, because as implied by Lemma~\ref{lem:combinatorics2}, the right-hand side grows as $(-1)^ke^{k(\epsilon-\frac{L_1}{2})}k^{-2}$ as $k$ increases when $L_1<2\epsilon$. 

By applying these expansion identities to \eqref{expansion of V_{KB}}, we obtain the proposition. 
\end{proof}

\section{Towards a construction of a recursive formula}\label{sec:higherTop}

In this section, we first discuss the difficulties in geometrically deriving a recursive formula for general topologies. These issues originate from the complexity of Gendulphe's regularisation and, in particular, do not arise in the orientable setting. We then shift our perspective from hyperbolic geometry to refined topological recursion, which may provide a new approach to deriving a recursive formula for Gendulphe--Norbury volumes.


\subsection{Subtleties for general topologies}

Similar to Mirzakhani's approach for $V_{g,n}^{+}(L_{[n]})$, one may attempt to derive a recursive formula for $V_{g,n}^{\epsilon,-}(L_{[n]})$ by applying Norbury--McShane identities. Following the procedure summarised in Section~\ref{sec:unfolding}, one can verify that \hyperref[item:Step1]{Step 1} and \hyperref[item:Step2]{Step 2} carry over unchanged, as the Norbury--McShane identities hold regardless of the regularisation. The main obstruction arises in \hyperref[item:Step3]{Step 3}: the integration over the cover $\mathcal{M}^{\gamma;\,\epsilon,-}_{g,n}(L_{[n]})$ can no longer be decomposed into a sum of products of volumes of lower-complexity moduli spaces times contributions from either a pair of pants or a two-bordered real projective plane.

To illustrate the issue, consider Figure~\ref{fig:1/2,3}, which depicts a possible contribution to $V_{\frac12,3}^{\epsilon,-}(L_{[3]})$.

\begin{figure}[htbp]
    \centering
    \begin{tikzpicture}[ 
        scale=0.3, 
        transform shape, 
        nodes={scale=2}, 
        thick, 
        line cap=round, 
        line join=round 
    ] 
        
        \coordinate (MiddleLeftTop) at (-4.07,5.05); 
        \coordinate (MiddleLeftBottom) at (-4.17,2.28); 
        \coordinate (MiddleLeftFront) at (-3.80,2.98); 
        \coordinate (MiddleLeftBack) at (-4.47,4.00); 
        \coordinate (MiddleRightTop) at (4.17,5.02); 
        \coordinate (MiddleRightBottom) at (4.07,2.22); 
        \coordinate (MiddleRightFront) at (3.78,3.02); 
        \coordinate (MiddleRightBack) at (4.47,4.00); 
        \coordinate (MiddleCrossLeft) at (-1.70506,6.88839); 
        \coordinate (MiddleCrossRight) at (2.34506,7.17161); 
        \coordinate (MiddleCrossTop) at (0.26768,7.77817); 
        \coordinate (MiddleCrossBottom) at (0.37232,6.28183); 
        \coordinate (MiddleChordOneLeft) at (-1.26491,7.40244); 
        \coordinate (MiddleChordOneRight) at (1.90491,6.65756); 
        \coordinate (MiddleChordTwoLeft) at (-1.19765,6.44061); 
        \coordinate (MiddleChordTwoRight) at (1.83765,7.61939); 
        
        \draw[red,thick,dashed] (MiddleCrossTop) .. controls (0.57,7.00) and (1.50,5.82) .. (2.35,5.28) .. controls (3.07,4.78) and (3.77,4.30) .. (MiddleRightBack); 
        
        \draw[red,thick] (MiddleCrossBottom) .. controls (1.20, 4.80) and (2.50, 3.20) .. (MiddleRightFront);
        
        \draw (MiddleLeftTop) .. controls (-3.27,4.97) and (-2.43,5.72) .. (MiddleCrossLeft); 
        \draw (MiddleCrossRight) .. controls (2.49,6.22) and (3.33,5.22) .. (MiddleRightTop); 
        \draw (MiddleLeftBottom) .. controls (-1.83,2.14) and (-0.93,2.51) .. (0.27,2.42) .. controls (1.73,2.31) and (2.97,2.13) .. (MiddleRightBottom); 
        \draw (0.32,7.03) ellipse[ x radius=2.03, y radius=0.75, rotate=4 ]; 
        \draw (MiddleChordOneLeft) -- (MiddleChordOneRight); 
        \draw (MiddleChordTwoLeft) -- (MiddleChordTwoRight); 
        \draw (MiddleLeftTop) .. controls (-4.58,4.56) and (-4.69,2.78) .. (MiddleLeftBottom); 
        \draw (MiddleLeftBottom) .. controls (-3.65,2.78) and (-3.51,4.55) .. (MiddleLeftTop); 
        \draw (MiddleRightTop) .. controls (4.67,4.54) and (4.58,2.76) .. (MiddleRightBottom); 
        \draw (MiddleRightBottom) .. controls (3.57,2.78) and (3.64,4.55) .. (MiddleRightTop); 
        \draw[black,thick] (0.32,7.03) ellipse[ x radius=2.03, y radius=0.75, rotate=4 ]; 
        \draw[black,thick] (MiddleLeftBottom) .. controls (-3.65,2.78) and (-3.51,4.55) .. (MiddleLeftTop); 
        \draw[black,thick] (MiddleRightBottom) .. controls (3.57,2.78) and (3.64,4.55) .. (MiddleRightTop); 
        \draw[black,thick] (MiddleLeftBottom) .. controls (-1.83,2.14) and (-0.93,2.51) .. (0.27,2.42) .. controls (1.73,2.31) and (2.97,2.13) .. (MiddleRightBottom); 
        
        \coordinate (SecondLeftTop) at (5.93,5.02); 
        \coordinate (SecondLeftBottom) at (5.83,2.22); 
        \coordinate (SecondSolidEnd) at (6.23,3.02); 
        \coordinate (SecondDashedEnd) at (5.54,4.00); 
        \coordinate (SecondRedStart) at (9.32345,3.34275); 
        
        \draw (10.29813,1.92836) .. controls (10.61631,1.54917) and (9.08422,0.26360) .. (8.76604,0.64279); 
        \draw (10.29813,1.92836) .. controls (9.97995,2.30755) and (8.44786,1.02198) .. (8.76604,0.64279); 
        \draw (9.21295,4.77736) .. controls (9.58580,5.26645) and (7.87221,6.53077) .. (7.49936,6.04168); 
        \draw (9.21295,4.77736) .. controls (8.84010,4.28827) and (7.12651,5.55259) .. (7.49936,6.04168); 
        \draw (SecondLeftTop) .. controls (6.43,4.54) and (6.34,2.76) .. (SecondLeftBottom); 
        \draw (SecondLeftBottom) .. controls (5.33,2.78) and (5.40,4.55) .. (SecondLeftTop); 
        \draw (8.76604,0.64279) .. controls (7.88,1.65) and (6.68,2.12) .. (SecondLeftBottom); 
        \draw (7.49936,6.04168) .. controls (6.92,5.34) and (6.32,5.12) .. (SecondLeftTop); 
        \draw[red,thick,dashed] (SecondRedStart) .. controls (9.58,3.48) and (7.18,4.08) .. (SecondDashedEnd); 
        \draw[red,thick] (SecondRedStart) .. controls (8.28,2.80) and (6.92,2.77) .. (SecondSolidEnd); 
        \draw[black,thick] (10.29813,1.92836) .. controls (9.43037,2.96252) and (8.92847,3.71624) .. (9.21295,4.77736); 
        \draw[black,thick] (SecondLeftBottom) .. controls (5.33,2.78) and (5.40,4.55) .. (SecondLeftTop); 
        
        
        \node at (-5.4, 3.6) {\Large $L_1$}; 
        
        \node at (8.3, 7.3) {\Large $L_2$};  
        \node at (9.5, -0.6) {\Large $L_3$}; 
        
        \node at (4.95, 1.2) {\Large $(\ell_1, \theta_1)$};  
        
        \node at (0.32, 8.8) {\Large $a$}; 
        
    \end{tikzpicture}
    \caption{A type of contribution to $V_{\frac12,3}^{\epsilon,-}(L_{[3]})$}\label{fig:1/2,3}
\end{figure}
\noindent Here $(\ell_1,\theta_1)$ are the Fenchel--Nielsen coordinates associated with the gluing, and $a$ is one of the 1-sided curves whose length is required by definition to be at least $\epsilon$. Observe that the two red arcs in each component are not closed curves, hence no regularisation is imposed on them, and their lengths can become arbitrarily small in a certain domain of $(\ell_1, a)$, prior to gluing. After gluing (possibly with a twist $\theta_1$), however, these two arcs combine to form a new 1-sided simple closed curve. To ensure that this curve also has length at least $\epsilon$, one must impose a constraint on $(\ell_1,\theta_1)$ which depends on $L_1,L_2,L_3$ and $a$.

We now turn to the next step in the recursion by considering surfaces of Euler characteristic $2-2g-n=-3$. As an example, Figure~\ref{fig:1/2,4-1} illustrates one possible contribution to $V_{\frac12,4}^{\epsilon,-}(L_{[4]})$.

\begin{figure}[htbp]
    \centering
    \begin{tikzpicture}[ 
        scale=0.3, 
        transform shape, 
        nodes={scale=2}, 
        thick, 
        line cap=round, 
        line join=round 
    ] 
        \coordinate (FirstRightTop) at (-5.83,5.05); 
        \coordinate (FirstRightBottom) at (-5.93,2.28); 
        \coordinate (FirstSolidEnd) at (-6.23,2.98); 
        \coordinate (FirstDashedEnd) at (-5.56,4.00); 
        \coordinate (FirstRedStart) at (-9.32345,3.34275); 
        
        \draw (-10.29813,1.92836) .. controls (-10.61631,1.54917) and (-9.08422,0.26360) .. (-8.76604,0.64279); 
        \draw (-10.29813,1.92836) .. controls (-9.97995,2.30755) and (-8.44786,1.02198) .. (-8.76604,0.64279); 
        \draw (-9.21295,4.77736) .. controls (-9.58580,5.26645) and (-7.87221,6.53077) .. (-7.49936,6.04168); 
        \draw (-9.21295,4.77736) .. controls (-8.84010,4.28827) and (-7.12651,5.55259) .. (-7.49936,6.04168); 
        \draw (FirstRightTop) .. controls (-6.34,4.56) and (-6.45,2.78) .. (FirstRightBottom); 
        \draw (FirstRightBottom) .. controls (-5.41,2.78) and (-5.27,4.55) .. (FirstRightTop); 
        \draw (-8.76604,0.64279) .. controls (-7.88,1.65) and (-6.68,2.15) .. (FirstRightBottom); 
        \draw (-7.49936,6.04168) .. controls (-6.92,5.36) and (-6.32,5.14) .. (FirstRightTop); 
        \draw[red,thick,dashed] (FirstRedStart) .. controls (-9.58,3.48) and (-7.18,4.08) .. (FirstDashedEnd); 
        \draw[red,thick] (FirstRedStart) .. controls (-8.28,2.80) and (-6.92,2.73) .. (FirstSolidEnd); 
        \draw[black,thick] (-10.29813,1.92836) .. controls (-9.43037,2.96252) and (-8.92847,3.71624) .. (-9.21295,4.77736); 
        \draw[black,thick] (FirstRightBottom) .. controls (-5.41,2.78) and (-5.27,4.55) .. (FirstRightTop); 
        
        \coordinate (MiddleLeftTop) at (-4.07,5.05); 
        \coordinate (MiddleLeftBottom) at (-4.17,2.28); 
        \coordinate (MiddleLeftFront) at (-3.80,2.98); 
        \coordinate (MiddleLeftBack) at (-4.47,4.00); 
        \coordinate (MiddleRightTop) at (4.17,5.02); 
        \coordinate (MiddleRightBottom) at (4.07,2.22); 
        \coordinate (MiddleRightFront) at (3.78,3.02); 
        \coordinate (MiddleRightBack) at (4.47,4.00); 
        \coordinate (MiddleCrossLeft) at (-1.70506,6.88839); 
        \coordinate (MiddleCrossRight) at (2.34506,7.17161); 
        \coordinate (MiddleCrossTop) at (0.26768,7.77817); 
        \coordinate (MiddleCrossBottom) at (0.37232,6.28183); 
        \coordinate (MiddleChordOneLeft) at (-1.26491,7.40244); 
        \coordinate (MiddleChordOneRight) at (1.90491,6.65756); 
        \coordinate (MiddleChordTwoLeft) at (-1.19765,6.44061); 
        \coordinate (MiddleChordTwoRight) at (1.83765,7.61939); 
        \coordinate (MiddleLowerTransition) at (0.05,2.42); 
        
        \draw[red,thick,dashed] (MiddleCrossTop) .. controls (-1.0, 5.5) and (-2.5, 4.00) .. (MiddleLeftBack); 
        
        \draw[red,thick] (MiddleLeftFront) .. controls (-2.73,3.26) and (-1.90,3.60) .. (-1.23,4.26) .. controls (-0.60,4.89) and (0.03,5.72) .. (MiddleCrossBottom); 
        
        \draw (MiddleLeftTop) .. controls (-3.27,4.97) and (-2.43,5.72) .. (MiddleCrossLeft); 
        \draw (MiddleCrossRight) .. controls (2.49,6.22) and (3.33,5.22) .. (MiddleRightTop); 
        \draw (MiddleLeftBottom) .. controls (-1.83,2.14) and (-0.93,2.51) .. (0.27,2.42) .. controls (1.73,2.31) and (2.97,2.13) .. (MiddleRightBottom); 
        \draw (0.32,7.03) ellipse[ x radius=2.03, y radius=0.75, rotate=4 ]; 
        \draw (MiddleChordOneLeft) -- (MiddleChordOneRight); 
        \draw (MiddleChordTwoLeft) -- (MiddleChordTwoRight); 
        \draw (MiddleLeftTop) .. controls (-4.58,4.56) and (-4.69,2.78) .. (MiddleLeftBottom); 
        \draw (MiddleLeftBottom) .. controls (-3.65,2.78) and (-3.51,4.55) .. (MiddleLeftTop); 
        \draw (MiddleRightTop) .. controls (4.67,4.54) and (4.58,2.76) .. (MiddleRightBottom); 
        \draw (MiddleRightBottom) .. controls (3.57,2.78) and (3.64,4.55) .. (MiddleRightTop); 
        \draw[black,thick] (0.32,7.03) ellipse[ x radius=2.03, y radius=0.75, rotate=4 ]; 
        \draw[black,thick] (MiddleLeftBottom) .. controls (-3.65,2.78) and (-3.51,4.55) .. (MiddleLeftTop); 
        \draw[black,thick] (MiddleRightBottom) .. controls (3.57,2.78) and (3.64,4.55) .. (MiddleRightTop); 
        \draw[black,thick] (MiddleLeftBottom) .. controls (-1.83,2.14) and (-0.93,2.51) .. (0.27,2.42) .. controls (1.73,2.31) and (2.97,2.13) .. (MiddleRightBottom); 
        
        \coordinate (SecondLeftTop) at (5.93,5.02); 
        \coordinate (SecondLeftBottom) at (5.83,2.22); 
        \coordinate (SecondSolidEnd) at (6.23,3.02); 
        \coordinate (SecondDashedEnd) at (5.54,4.00); 
        \coordinate (SecondRedStart) at (9.32345,3.34275); 
        
        \draw (10.29813,1.92836) .. controls (10.61631,1.54917) and (9.08422,0.26360) .. (8.76604,0.64279); 
        \draw (10.29813,1.92836) .. controls (9.97995,2.30755) and (8.44786,1.02198) .. (8.76604,0.64279); 
        \draw (9.21295,4.77736) .. controls (9.58580,5.26645) and (7.87221,6.53077) .. (7.49936,6.04168); 
        \draw (9.21295,4.77736) .. controls (8.84010,4.28827) and (7.12651,5.55259) .. (7.49936,6.04168); 
        \draw (SecondLeftTop) .. controls (6.43,4.54) and (6.34,2.76) .. (SecondLeftBottom); 
        \draw (SecondLeftBottom) .. controls (5.33,2.78) and (5.40,4.55) .. (SecondLeftTop); 
        \draw (8.76604,0.64279) .. controls (7.88,1.65) and (6.68,2.12) .. (SecondLeftBottom); 
        \draw (7.49936,6.04168) .. controls (6.92,5.34) and (6.32,5.12) .. (SecondLeftTop); 
        \draw[black,thick] (10.29813,1.92836) .. controls (9.43037,2.96252) and (8.92847,3.71624) .. (9.21295,4.77736); 
        \draw[black,thick] (SecondLeftBottom) .. controls (5.33,2.78) and (5.40,4.55) .. (SecondLeftTop); 
        
        \node at (-8.3, 7.3) {\Large $L_1$}; 
        \node at (-9.5, -0.6) {\Large $L_4$}; 
        \node at (8.3, 7.3) {\Large $L_2$};  
        \node at (9.5, -0.6) {\Large $L_3$};  
        
        \node at (-5.05, 1.2) {\Large $(\ell_2, \theta_2)$}; 
        \node at (4.95, 1.2) {\Large $(\ell_1, \theta_1)$};  
        
        \node at (0.32, 8.8) {\Large $a$}; 
    \end{tikzpicture}\caption{Some contributions to $V_{\frac12,4}^{\epsilon,-}(L_{[4]})$}\label{fig:1/2,4-1}
\end{figure}

\noindent As before, $(\ell_1,\theta_1)$ must satisfy the constraint discussed above, with $L_1$ replaced by $\ell_2$, to ensure that the corresponding 1-sided curve has length greater than $\epsilon$. Additionally, a new constraint must be imposed on $(\ell_2,\theta_2)$ to regularise the red curve in Figure~\ref{fig:1/2,4-1}. This constraint will depend not only on $L_1$ and $L_4$, but also on $\ell_1$ and $a$. Consequently, the two constraints are mutually dependent, so the admissible domain for $(\ell_1,\theta_1,a)$ differs from that in the case $(g,n)=(\frac12,3)$. This shows that \hyperref[item:Step3]{Step 3} does not extend straightforwardly to the non-orientable setting.

Furthermore, a new 1-sided curve can also take the following form:

\begin{figure}[htbp]
   \centering
   \begin{tikzpicture}[ 
       scale=0.3, 
       transform shape, 
       nodes={scale=2}, 
       thick, 
       line cap=round, 
       line join=round 
   ] 
       \coordinate (FirstRightTop) at (-5.83,5.05); 
       \coordinate (FirstRightBottom) at (-5.93,2.28); 
       \coordinate (FirstSolidEnd) at (-6.23,2.98); 
       \coordinate (FirstDashedEnd) at (-5.56,4.00); 
       \coordinate (FirstRedStart) at (-9.32345,3.34275); 
       
       \draw (-10.29813,1.92836) .. controls (-10.61631,1.54917) and (-9.08422,0.26360) .. (-8.76604,0.64279); 
       \draw (-10.29813,1.92836) .. controls (-9.97995,2.30755) and (-8.44786,1.02198) .. (-8.76604,0.64279); 
       \draw (-9.21295,4.77736) .. controls (-9.58580,5.26645) and (-7.87221,6.53077) .. (-7.49936,6.04168); 
       \draw (-9.21295,4.77736) .. controls (-8.84010,4.28827) and (-7.12651,5.55259) .. (-7.49936,6.04168); 
       \draw (FirstRightTop) .. controls (-6.34,4.56) and (-6.45,2.78) .. (FirstRightBottom); 
       \draw (FirstRightBottom) .. controls (-5.41,2.78) and (-5.27,4.55) .. (FirstRightTop); 
       \draw (-8.76604,0.64279) .. controls (-7.88,1.65) and (-6.68,2.15) .. (FirstRightBottom); 
       \draw (-7.49936,6.04168) .. controls (-6.92,5.36) and (-6.32,5.14) .. (FirstRightTop); 
       \draw[red,thick,dashed] (FirstRedStart) .. controls (-9.58,3.48) and (-7.18,4.08) .. (FirstDashedEnd); 
       \draw[red,thick] (FirstRedStart) .. controls (-8.28,2.80) and (-6.92,2.73) .. (FirstSolidEnd); 
       \draw[black,thick] (-10.29813,1.92836) .. controls (-9.43037,2.96252) and (-8.92847,3.71624) .. (-9.21295,4.77736); 
       \draw[black,thick] (FirstRightBottom) .. controls (-5.41,2.78) and (-5.27,4.55) .. (FirstRightTop); 
       
       \coordinate (MiddleLeftTop) at (-4.07,5.05); 
       \coordinate (MiddleLeftBottom) at (-4.17,2.28); 
       \coordinate (MiddleLeftFront) at (-3.80,2.98); 
       \coordinate (MiddleLeftBack) at (-4.47,4.00); 
       \coordinate (MiddleRightTop) at (4.17,5.02); 
       \coordinate (MiddleRightBottom) at (4.07,2.22); 
       \coordinate (MiddleRightFront) at (3.78,3.02); 
       \coordinate (MiddleRightBack) at (4.47,4.00); 
       \coordinate (MiddleCrossLeft) at (-1.70506,6.88839); 
       \coordinate (MiddleCrossRight) at (2.34506,7.17161); 
       \coordinate (MiddleCrossTop) at (0.26768,7.77817); 
       \coordinate (MiddleCrossBottom) at (0.37232,6.28183); 
       \coordinate (MiddleChordOneLeft) at (-1.26491,7.40244); 
       \coordinate (MiddleChordOneRight) at (1.90491,6.65756); 
       \coordinate (MiddleChordTwoLeft) at (-1.19765,6.44061); 
       \coordinate (MiddleChordTwoRight) at (1.83765,7.61939); 
       \coordinate (MiddleLowerTransition) at (0.05,2.42); 
       
       \draw[red,thick,dashed] (MiddleCrossTop) .. controls (0.57,7.00) and (1.50,5.82) .. (2.35,5.28) .. controls (3.07,4.78) and (3.77,4.30) .. (MiddleRightBack); 
       \draw[red,thick,dashed] (MiddleLowerTransition) .. controls (-0.83,2.51) and (-1.63,2.98) .. (-2.31,3.48) .. controls (-3.07,4.02) and (-3.83,4.00) .. (MiddleLeftBack); 
       \draw[red,thick] (MiddleLeftFront) .. controls (-2.73,3.26) and (-1.90,3.60) .. (-1.23,4.26) .. controls (-0.60,4.89) and (0.03,5.72) .. (MiddleCrossBottom); 
       \draw[red,thick] (MiddleRightFront) .. controls (2.97,3.07) and (2.39,3.19) .. (1.76,3.03) .. controls (1.09,2.85) and (0.57,2.54) .. (MiddleLowerTransition); 
       \draw (MiddleLeftTop) .. controls (-3.27,4.97) and (-2.43,5.72) .. (MiddleCrossLeft); 
       \draw (MiddleCrossRight) .. controls (2.49,6.22) and (3.33,5.22) .. (MiddleRightTop); 
       \draw (MiddleLeftBottom) .. controls (-1.83,2.14) and (-0.93,2.51) .. (0.27,2.42) .. controls (1.73,2.31) and (2.97,2.13) .. (MiddleRightBottom); 
       \draw (0.32,7.03) ellipse[ x radius=2.03, y radius=0.75, rotate=4 ]; 
       \draw (MiddleChordOneLeft) -- (MiddleChordOneRight); 
       \draw (MiddleChordTwoLeft) -- (MiddleChordTwoRight); 
       \draw (MiddleLeftTop) .. controls (-4.58,4.56) and (-4.69,2.78) .. (MiddleLeftBottom); 
       \draw (MiddleLeftBottom) .. controls (-3.65,2.78) and (-3.51,4.55) .. (MiddleLeftTop); 
       \draw (MiddleRightTop) .. controls (4.67,4.54) and (4.58,2.76) .. (MiddleRightBottom); 
       \draw (MiddleRightBottom) .. controls (3.57,2.78) and (3.64,4.55) .. (MiddleRightTop); 
       \draw[black,thick] (0.32,7.03) ellipse[ x radius=2.03, y radius=0.75, rotate=4 ]; 
       \draw[black,thick] (MiddleLeftBottom) .. controls (-3.65,2.78) and (-3.51,4.55) .. (MiddleLeftTop); 
       \draw[black,thick] (MiddleRightBottom) .. controls (3.57,2.78) and (3.64,4.55) .. (MiddleRightTop); 
       \draw[black,thick] (MiddleLeftBottom) .. controls (-1.83,2.14) and (-0.93,2.51) .. (0.27,2.42) .. controls (1.73,2.31) and (2.97,2.13) .. (MiddleRightBottom); 
       
       \coordinate (SecondLeftTop) at (5.93,5.02); 
       \coordinate (SecondLeftBottom) at (5.83,2.22); 
       \coordinate (SecondSolidEnd) at (6.23,3.02); 
       \coordinate (SecondDashedEnd) at (5.54,4.00); 
       \coordinate (SecondRedStart) at (9.32345,3.34275); 
       
       \draw (10.29813,1.92836) .. controls (10.61631,1.54917) and (9.08422,0.26360) .. (8.76604,0.64279); 
       \draw (10.29813,1.92836) .. controls (9.97995,2.30755) and (8.44786,1.02198) .. (8.76604,0.64279); 
       \draw (9.21295,4.77736) .. controls (9.58580,5.26645) and (7.87221,6.53077) .. (7.49936,6.04168); 
       \draw (9.21295,4.77736) .. controls (8.84010,4.28827) and (7.12651,5.55259) .. (7.49936,6.04168); 
       \draw (SecondLeftTop) .. controls (6.43,4.54) and (6.34,2.76) .. (SecondLeftBottom); 
       \draw (SecondLeftBottom) .. controls (5.33,2.78) and (5.40,4.55) .. (SecondLeftTop); 
       \draw (8.76604,0.64279) .. controls (7.88,1.65) and (6.68,2.12) .. (SecondLeftBottom); 
       \draw (7.49936,6.04168) .. controls (6.92,5.34) and (6.32,5.12) .. (SecondLeftTop); 
       \draw[red,thick,dashed] (SecondRedStart) .. controls (9.58,3.48) and (7.18,4.08) .. (SecondDashedEnd); 
       \draw[red,thick] (SecondRedStart) .. controls (8.28,2.80) and (6.92,2.77) .. (SecondSolidEnd); 
       \draw[black,thick] (10.29813,1.92836) .. controls (9.43037,2.96252) and (8.92847,3.71624) .. (9.21295,4.77736); 
       \draw[black,thick] (SecondLeftBottom) .. controls (5.33,2.78) and (5.40,4.55) .. (SecondLeftTop); 
       
       \node at (-8.3, 7.3) {\Large $L_1$}; 
       \node at (-9.5, -0.6) {\Large $L_4$}; 
       \node at (8.3, 7.3) {\Large $L_3$};  
       \node at (9.5, -0.6) {\Large $L_3$};  
       
       \node at (-5.05, 1.2) {\Large $(\ell_2, \theta_2)$}; 
       \node at (4.95, 1.2) {\Large $(\ell_1, \theta_1)$};  
       
       \node at (0.32, 8.8) {\Large $a$}; 
   \end{tikzpicture}\caption{Other contributions to $V_{\frac12,4}^{\epsilon,-}(L_{[4]})$}\label{fig:1/2,4-2}
\end{figure}

\noindent The new 1-sided curve in Figure~\ref{fig:1/2,4-2} propagates through all three components, thereby imposing an additional global constraint on  the Fenchel--Nielsen coordinates $(\ell_1,\theta_1,a,\ell_2,\theta_2)$.


The situation becomes increasingly complicated as $2-2g-n$ decreases; in particular, a different type of constraint must also be imposed when forming non-orientable handles. Moreover, similar to Figure \ref{fig:1/2,4-2}, the relevant constraints are no longer local in general: one must bound the lengths of $1$-sided curves that cross the gluing seam(s) and propagate through more than the two adjacent building pieces (of topology $(0,3)$ or $(\frac12,2)$) connected by the seam(s). Furthermore, when the gluing curves are sufficiently long, one must also account for the possibility that, depending on the twist coordinates, 1-sided curves crossing the gluing seams more than twice can become shorter than those crossing them only twice.




Combining these aspects, it appears highly challenging to geometrically derive an explicit recursion in the non-orientable setting.

\subsubsection{A recursive mechanism from hyperbolic geometry} 


Despite these obstacles, the geometry still suggests an underlying recursive mechanism.
Observe first that the only new 1-sided curves that can become shorter than~$\epsilon$ are those crossing the gluing curve, whose Fenchel--Nielsen coordinates are $(\ell,\theta)$. Recall that, since $\epsilon$ is assumed sufficiently small in the sense of Section~\ref{sec:reg}, no restriction needs to be imposed on the twist coordinate whenever the corresponding gluing seam satisfies $\ell<\epsilon$. In this region of the moduli space, \hyperref[item:Step3]{Step 3} can still be extended to the non-orientable setting without further restrictions, thereby yielding the corresponding contribution to the regularised volume. 

In the regime $\ell>\epsilon$, on the other hand, one must remove from the standard twist-coordinate domain $[0,\ell)$ all intervals that would allow the formation of new 1-sided geodesics of length smaller than $\epsilon$. Using the assumption $\sinh (\epsilon / 2) < 1$, one can show that, for each gluing length~$\ell$, any potentially short new 1-sided geodesic crosses the gluing seam only a bounded number of times. Therefore, only finitely many inadmissible intervals must be removed from the twist-coordinate domain. These inadmissible intervals depend not only on the geometry of the pair of pants being attached, but also on the Fenchel--Nielsen coordinates of the lower-complexity surfaces. Although this procedure is highly implicit and appears difficult to exploit in practice, a more detailed analysis may yield useful structural results for the Gendulphe--Norbury volumes.


In summary, while a recursive structure exists in the non-orientable setting, a naive application of Mirzakhani's technique based solely on the Norbury--McShane identities is unlikely to give rise to an explicit formulation. This motivates seeking a more tractable formulation via refined topological recursion, to which we now turn.

\subsection{Refined topological recursion}\label{sec:RTR}
For a compact Riemann surface $\Sigma$, refined topological recursion, originally proposed in \cite{CE06} and revisited in \cite{KO22,Osu23-1}, constructs a family of correlators $\omega_{g,n}$\footnote{For brevity of notation, we drop the $\mathfrak{b}$-dependence from $\omega_{g,n}$, in contrast to the notation $\omega_{g,n}^\mathfrak{b}$ used in Section~\ref{sec:intro}.}, namely meromorphic sections of $K_\Sigma^{\boxtimes n}$, indexed by $g\in\frac12\mathbb{Z}_{\geq0}$ and $n\in\mathbb{Z}_{\geq1}$. When $\Sigma=\mathbb{P}^1$, all expected properties have been proved in  \cite{KO22,Osu23-1,Osu23-2}, while some of them remain conjectural when $H^1(\Sigma)$ is nontrivial. For our purposes, we need to extend the framework to a non-compact case $\Sigma=\mathbb{C}$.
 
In the standard setting, the \emph{spectral curve} $\mathcal{S}_{\text{WP}}$ for Weil--Petersson volumes is given by the tuple $(\Sigma,x,y,\omega_{0,2})$, where $\Sigma=\mathbb{C}$\footnote{In the literature, $\Sigma$ is often chosen to be $\mathbb{P}^1$. There is no practical difference in the orientable setting, because $\omega_{g,n}$ are rational without poles at $\infty$, hence they can be viewed as meromorphic multidifferentials on $\mathbb{C}$ or $\mathbb{P}^1$. In the refined setting, however, the $\omega_{g,n}$ have an essential singularity at $\infty$, so it is in fact more natural to set $\Sigma=\mathbb{C}$.}, and the remaining objects are defined in terms of the standard global coordinate $z$ as follows:
\begin{equation}
x(z)=\frac{z^2}{2},\quad y(z)=-\frac{\sin2\pi z}{2\pi},\quad\omega_{0,2}(z_1,z_2)=\frac{dz_1dz_2}{(z_1+z_2)^2},
\end{equation}
where we adopt a slightly different convention for $\omega_{0,2}$ compared to the one used in \cite{EO08}. We denote by $\sigma:\Sigma\to\Sigma$ the involution operator whose action on $(x,y)$ is given by $\sigma:(x,y)\mapsto(x,-y)$. In terms of the global coordinate $z$, this simply becomes $\sigma(z)=-z$. We now upgrade $\mathcal{S}_{\text{WP}}$ to the refined setting:

\begin{definition}\label{def:curve}
Fix $\mathfrak{b}\in\mathbb{C}$ and $\epsilon\in\mathbb{R}_{>0}$. The \emph{refined spectral curve} $\mathcal{S}^{\epsilon}_{\textup{total}}$ for regularised volumes of moduli spaces of hyperbolic Klein surfaces is a tuple $(\Sigma,x,y,\omega_{0,2},\omega_{\frac12,1})$ such that $\Sigma,x,y,\omega_{0,2}$ are the same as above, and $\omega_{\frac12,1}$ is set to:
\begin{equation}
\omega_{\frac12,1}(z)\coloneqq\frac{\mathfrak{b}}{2}\left(-\frac{dy(z)}{y(z)}+\frac{e^{2z\epsilon}-e^{-2z\epsilon}}{2z}dz+\sum_{k\geq1}\left(\frac{e^{2(z-\frac{k}{2})\epsilon}}{z-\frac{k}{2}}-\frac{e^{-2(z+\frac{k}{2})\epsilon}}{z+\frac{k}{2}}\right)dz\right)\label{w_{1/2,1}}.
\end{equation}
\end{definition}

In Section \ref{sec:evidence}, we show that the refined spectral curve above indeed reproduces the correct volumes of the moduli spaces of two-bordered real projective planes and one-bordered Klein bottles. The expression \eqref{w_{1/2,1}} was discovered through educated experimentation by the authors, and it remains an interesting open question whether it can be \emph{derived} from the perspective of hyperbolic geometry.

\begin{remark}\label{rem:RSC}
Let us justify that the definition above sits inside a natural extension of the framework of \cite{Osu23-1}. When $\Sigma$ is a compact Riemann surface, it is required in \cite{Osu23-1} that $\omega_{\frac12,1}$ admits the properties below:
\begin{itemize}
\item[{RSC1}:] \label{item:RSC1} it can have residues at zeroes and poles of $ydx$ but nowhere else, and no higher order poles,
\item[{RSC2}:] \label{item:RSC2} one can add to it holomorphic differentials when $H^1(\Sigma)$ is nontrivial,
\item[{RSC3}:] \label{item:RSC3}$\omega_{\frac12,1}+\frac{\mathfrak{b}}{2}\frac{dy}{y}$ is anti-invariant under the involution $\sigma$.
\end{itemize}
It is evident that $\omega_{\frac12,1}$ in \eqref{w_{1/2,1}} satisfies the first and third conditions; in particular, for $k\geq1$:
\begin{equation*}
\Res_{z=\frac{k}{2}}\,\omega_{\frac12,1}(z)=0,\quad \Res_{z=-\frac{k}{2}}\,\omega_{\frac12,1}(z)=-\mathfrak{b},\quad \Res_{z=0}\,\omega_{\frac12,1}(z)=-\frac{\mathfrak{b}}{2}.
\end{equation*}
Moreover, in contrast to the case $\Sigma=\mathbb{P}^1$, the point $z=\infty$ does not belong to $\Sigma$. Consequently, the second property permits the exponential factors $e^{\pm2(z\mp\frac{k}{2})\epsilon}$ to appear in the definition. Last but not least, the sum in $k$ is absolutely convergent precisely because $\epsilon\in\mathbb{R}_{>0}$. 
\end{remark}

Extending the approach of \cite{Osu23-1}, we define refined topological recursion as below.

\begin{definition}\label{def:RTR}
Given $\mathcal{S}^{\epsilon}_{\textup{total}}$, the \emph{refined topological recursion} recursively constructs the stable correlators $\omega_{g,n+1}$, which are meromorphic sections of $K_\Sigma^{\boxtimes (n+1)}$, for $g\in\frac12\mathbb{Z}_{\geq0}$ and $n\in\mathbb{Z}_{\geq0}$ satisfying $2g-2+n+1\geq1$, via the formula
\begin{align}
\omega_{g,n+1}(z_0,z_{[n]})\coloneqq&\left(\sum_{k\geq1}\Res_{z=\frac{k}{2}}+\sum_{i=0}^n\Res_{z=z_i}-\Res_{z=0}-\sum_{k\geq1}\Res_{z=-\frac{k}{2}}-\sum_{i=0}^n\Res_{z=-z_i}\right)\frac{\eta^z(z_1)}{4\omega_{0,1}(z)}{\rm Rec}_{g,n+1}(z,z_{[n]}),\label{RTR}
\end{align}
where $\eta^z(z_1)\coloneqq\int_{-z}^z\omega_{0,2}(z_1,\cdot)$, and
\begin{align}
{\rm Rec}_{g,n+1}(z,z_{[n]})=&\sum_{\substack{g_1+g_2=g\\J_1\sqcup J_2=[n]}}^*\omega_{g_1,1+|J_1|}(z,z_{J_1})\,\omega_{g_2,1+|J_2|}(z,z_{J_2})+\sum_{i=1}^n\frac{dx(z)dx(z_i)}{(x(z)-x(z_i))^2}\omega_{g,n}(z,z_{[n]\backslash\{i\}})\nonumber\\
&+\omega_{g-1,n+2}(z,z,z_{[n]})+\mathfrak{b}\,dx(z)d_z\frac{\omega_{g-\frac12,n+1}(z,z_{[n]})}{dx(z)}.\label{Rec}
\end{align}
Here, the symbol $*$ in the summation indicates that terms involving $\omega_{0,1}$ are omitted, whereas those involving $\omega_{0,2}$ and $\omega_{\frac12,1}$ are retained. We also write $d_z$ for the exterior derivative with respect to $z$.
\end{definition}

Viewed as a function of $\mathfrak{b}$, one can show by induction that, for each $g\in\frac12\mathbb{Z}_{\geq0}$, $\omega_{g,n}$ is a polynomial in $\mathfrak{b}$ of degree at most $2g$, for every $n\in\mathbb{Z}_{\geq1}$. Moreover, when $\mathfrak{b}=0$, an induction argument analogous to that of \cite[Lemma 3.4]{Osu23-1} shows that $\omega_{g,n}$, as defined in Definition \ref{def:RTR}, coincides with the corresponding correlator from the Chekhov--Eynard--Orantin topological recursion \cite{CE05,EO07} for degree-two spectral curves.

However, a subtlety of Definition \ref{def:RTR} is that, since $\Sigma$ is non-compact, other properties proven in \cite{KO22,Osu23-1,Osu23-2} need to be re-checked. In addition, since the formula involves residues at infinitely many points, convergence issues must also be addressed. A complete treatment requires a careful analysis and lies beyond the scope of the present paper, so we leave it to future work. Let us instead summarise the expected properties of $\omega_{g,n}$, which we verify in the cases with $2g-2+n=1$ in Section \ref{sec:evidence}:
\begin{conjecture}\label{conj:RTR}
   For every $g\in\frac12\mathbb{Z}_{\geq0}$ and $n\in\mathbb{Z}_{\geq1}$ satisfying $2g-2+n\geq1$, let $\omega_{g,n}$ be the multidifferential defined in {\rm Definition \ref{def:RTR}}. Then $\omega_{g,n}$ satisfies the following properties:
    \begin{description}
     \item[{\rm RTR0}] \label{item:RTR0} the termwise inverse Laplace transform \textup{(Definition \ref{def:termwise})} of $\omega_{g,n}$ is convergent,
       \item[{\rm RTR1}] \label{item:RTR1}$\omega_{g,n}$ is a symmetric multidifferential viewed as a (possibly divergent) series,
       \item[{\rm RTR2}] \label{item:RTR2}$\omega_{g,n}$ has no residues with respect to any variable,
       \item[{\rm RTR3}] \label{item:RTR3}$\omega_{g,n}$ can have poles only at $z_i=-z_j$ and at $z_i=-\frac{k}{2}$ for all $i,j\in[n]$ and $k\in\mathbb{Z}_{\geq1}$.
   \end{description}
\end{conjecture}

Before turning to the relation between $V_{g,n}^\epsilon$ and $\omega_{g,n}$, we introduce some notation that will be convenient in what follows. For a multidifferential $\omega$, we define:
\begin{equation*}
\Delta_{z_i}\omega(z_1,...,z_n)\coloneqq\omega(z_1,...,z_i,...,z_n)-\omega(z_1,...,-z_i,...,z_n).
\end{equation*}
In particular, in our setting we have
\begin{align}
\Delta_z\omega_{0,2}(z,z_1)=&\;2\omega_{0,2}(z,z_1)+\frac{dx(z)dx(z_1)}{(x(z)-x(z_1))^2},\label{Delta_0,2}\\
\Delta_z\omega_{\frac12,1}(z)=&\;\mathfrak{b}\sum_{k\geq0}\left(\frac{e^{-2(z+\frac{k}{2})\epsilon}}{z+\frac{k}{2}}-\frac{e^{2(z-\frac{k}{2})\epsilon}}{z-\frac{k}{2}}\right)dz.
\end{align}

\subsection{Volumes via refined topological recursion}


When $\mathfrak{b}=0$, it was shown in \cite{EO07-2} that $\omega_{g,n}$ and $V_{g,n}^{+}$ are related by the Laplace transform $\mathcal{L}$. More concretely:
\begin{theorem}[\cite{EO07-2}]\label{thm:EOMir}
    Set $\mathfrak{b}=0$. For every $g\in\mathbb{Z}_{\geq0}$ and $n\in\mathbb{Z}_{\geq1}$ with $2g-2+n>0$, 
    \begin{equation}
     \prod_{i=1}^nL_i  \cdot V_{g,n}^{+}(L_{[n]})=(\mathcal{L}_{L_1}^{-1}\cdots \mathcal{L}_{L_n}^{-1}).\; \omega_{g,n}.
    \end{equation}
\end{theorem}

One may speculate that the total volumes $V^\epsilon_{g,n}$ are likewise related to the refined correlators $\omega_{g,n}$ by a Laplace transform for a certain choice of the refinement parameter $\mathfrak{b}$. Since the $\omega_{g,n}$ are possibly divergent series, however, the meaning of the inverse Laplace transform must be treated with care. To facilitate explicit computations, we first introduce the notion of a termwise Laplace transform.
\begin{definition}\label{def:termwise}
    Let $\omega$ be a series $\omega(z)=\sum_{k\geq0} \omega_{k}(z)$, where each $\omega_{k}(z)$ is a meromorphic differential on $\mathbb{C}$. For $L\in\mathbb{R}_{\geq0}$, we define the \emph{termwise inverse Laplace transform} $\hat{\mathcal{L}}^{-1}_L$ by
    \begin{equation*}
    \hat{\mathcal{L}}^{-1}_L.\;\omega\coloneqq\sum_{k\geq0} \mathcal{L}^{-1}_L.\,\omega_{k},
\end{equation*}
where $\mathcal{L}^{-1}_L$ denotes the usual inverse Laplace transform.
\end{definition}

If $\omega$ is a convergent series, it follows that ${\mathcal{L}}^{-1}_L.\;\omega=\hat{\mathcal{L}}^{-1}_L.\;\omega$. A crucial difference appears when $\omega$ is a divergent series: the termwise inverse Laplace transform $\hat{\mathcal{L}}^{-1}_L.\;\omega$ makes the divergent behaviour factorially better, and it may return a convergent series. This is related to \hyperref[item:RTR0]{RTR0} in Conjecture~\ref{conj:RTR}.

We now pose the following question as a step towards a generalisation of Theorem~\ref{thm:EOMir}:

\begin{question}\label{conj:main}
   Let $g\in\mathbb{Z}_{\geq0}$ and $n\in\mathbb{Z}_{\geq1}$ satisfy $2g-2+n>0$, and assume that Conjecture~\ref{conj:RTR} holds. Let $\omega_{g,n}$ denote the multidifferential produced by refined topological recursion on the refined spectral curve $\mathcal{S}^{\epsilon}_{\textup{total}}$, and let $V_{g,n}^\epsilon(L_{[n]})$ denote the total volume of the regularised moduli space of $n$-bordered Klein surfaces of genus $g$ with boundary lengths $L_{[n]}$.
 Introducing a new parameter $b$ by $\mathfrak{b} = -\frac{b}{\sqrt{1+b}}$, does there exist a relation of the following form between $V_{g,n}^\epsilon(L_{[n]})$ and $\omega_{g,n}$ when $b=1$?
    \begin{equation}
     \prod_{i=1}^nL_i\cdot V_{g,n}^\epsilon(L_{[n]})= (1+b)^{g}(\hat{\mathcal{L}}_{L_1}^{-1}\cdots \hat{\mathcal{L}}_{L_n}^{-1}).\;\omega_{g,n}\Big|_{b=1}.\label{V and omega relation}
    \end{equation}
\end{question}
The reparameterisation to $b$ and the normalisation $(1+b)^g$ are inspired by the relation between refined topological recursion and the so-called $b$-Hurwitz numbers \cite{CDO24-1,CDO24-2}. The parameter $b$ is related to the notion of measure of non-orientability due to Chapuy and Do\l\k{e}ga \cite{CD20}. It is suggested by these papers that when $b=1$ ($\mathfrak{b}=-\frac{1}{\sqrt{2}}$), refined topological recursion generates enumerative invariants that equally count contributions from orientable and non-orientable surfaces, which is analogous to our total volume $V_{g,n}^\epsilon$.

\begin{remark}
We will show in Section \ref{sec:evidence} that Question \ref{conj:main} has an affirmative answer when $2g-2+n=1$. From the geometric perspective, the corresponding volumes $V^\epsilon_{g,n}$ with $2g-2+n=1$ cannot be obtained recursively and must instead be computed directly. By contrast, the refined correlators $\omega_{g,n}$ with $2g-2+n=1$ are determined recursively by Definition~\ref{def:RTR} from the $\omega_{g,n}$ with $2g-2+n\leq0$. Thus, refined topological recursion has already passed at least one nontrivial check, providing some evidence that the total volumes $V_{g,n}^\epsilon(L_{[n]})$ may likewise admit a recursive description through the dictionary~\eqref{V and omega relation}.
\end{remark}

\subsubsection{Comments on poles along the anti-diagonal}

Due to the complicated pole structure of $\omega_{g,n}$, there exists an extra layer of subtleties on the action of $\hat{\mathcal{L}}^{-1}_L$ in the refined setting, which originates from poles along the anti-diagonal. Let us demonstrate it with a toy example. Take
\begin{equation*}
   \omega(z_1,z_2)=\frac{\left(z_1^4+3 z_1^3z_2+3 z_1^2 z_2^2+3 z_1z_2^3+z_2^4\right)}{2 z_1^3 z_2^3
   \left(z_1+z_2\right)^3}.
\end{equation*}
It is a symmetric bidifferential without residues, and in terms of $z_1$, it has poles at $z_1=0$ and $z_1=-z_2$. Then, one finds:
\begin{equation*}
\hat{\mathcal{L}}^{-1}_{L_2}.\;  \omega(z_1,z_2)= \frac{L_2^2}{4 z_1^2}+\frac{L_2^2 e^{-L_2 z_1}}{4 z_1^2}+\frac{L_2e^{-L_2 z_1}}{2 z_1^3}
\end{equation*}
and
\begin{equation}
  (\hat{\mathcal{L}}^{-1}_{L_1}\hat{\mathcal{L}}^{-1}_{L_2}).\; \omega(z_1,z_2)=\frac{L_2^2L_1}{4}+\frac{L_1 L_2(L_1-L_2)\theta(L_1-L_2)}{4}=\frac{L_1L_2\max(L_1,L_2)}{4},\label{L2L1}
\end{equation}
where $\theta(x)$ is the Heaviside step function. One can see that the step function appears as a consequence of the pole at $z_1+z_2=0$. 

Since the refined correlators $\omega_{g,n}$ have poles along the anti-diagonal, phenomena
similar to the toy example above naturally arise. In practice, however, it is sometimes convenient to impose some inequalities among the boundary lengths $L_{[n]}$ from the beginning, thereby avoiding the appearance of step functions. In particular, these inequalities allow one to choose the contours for the inverse Laplace transform uniformly, whereas otherwise the contours may need to be chosen separately for each term. A similar piecewise structure is discussed in \cite{CGGO26}.

\subsubsection{Comments on generic values of $b$}

Another motivation for introducing refined topological recursion in the present paper is that it naturally leads to the following intriguing question: what happens if the parameter $b$ is left unfixed? Motivated by this question, we \emph{define} $V_{g,n}^{\epsilon,b}(L_{[n]})$ by
\begin{equation}
     \prod_{i=1}^nL_i\cdot V_{g,n}^{\epsilon,b}(L_{[n]})\coloneqq (1+b)^{g}(\hat{\mathcal{L}}_{L_1}^{-1}\cdots \hat{\mathcal{L}}_{L_n}^{-1}).\;\omega_{g,n}.
    \end{equation}
Viewed as a function of $b$, each $V_{g,n}^{\epsilon,b}(L_{[n]})$ turns out to be a polynomial in $b$ of degree $2g$. For example, as shown in Section~\ref{sec:RTR Klein bottle}, we have
\[V_{1,1}^{\epsilon,b}(L_1)=\frac{L_1^2}{48}+\frac{\pi^2}{12}+b\left(\frac{L_1^2}{48}+\frac{\pi^2}{12}\right)+b^2\left(-\frac{L_1^2}{48}-\frac{\pi^2}{12}-{\rm Li}_2\left(-\frac{\cosh^2\frac{L_1}{4}}{\sinh^2\frac{\epsilon}{2}}\right)\right).\]

Using the identity \eqref{dilog identity}, one can show that the coefficient of $b^2$ is actually positive whenever $L_1>0$ and $0<\sinh\frac{\epsilon}{2}<1$. 
Moreover,
$$
V_{1,1}^{\epsilon,0}(L_1)=V_{1,1}^{+}(L_1),
$$ 
while
$$
V_{1,1}^{\epsilon,1}(L_1)=V_{1,1}^{\epsilon}(L_1)=V_{1,1}^{+}(L_1)+V_{1,1}^{\epsilon,-}(L_1)
$$
by \eqref{V_{0,3} and V_{1,1}} and Theorem~\ref{thm:KB}. Thus, $V_{1,1}^{\epsilon,b}(L_1)$ interpolates between the oriented and total volumes.  

This naturally raises the question of whether the parameter $b$ admits a geometric interpretation, analogous to the one given by the measure of non-orientability introduced by Chapuy and Do\l\k{e}ga in a more combinatorial context \cite{CD20}. It is worth mentioning that, in the limit $L_k\to\infty$, such a geometric interpretation is indeed available: the leading-order terms are (strictly speaking, conjecturally) realised as volumes of moduli spaces of metric Möbius graphs equipped with a measure of non-orientability \cite{CGGO26}.

\subsubsection{Comments on Stanford's proposal}

In the context of Jackiw--Teitelboim gravity in physics, Stanford \cite{Sta23} proposed a recursion formula for $\mathcal{V}_{g,n}^{\epsilon}(L_{[n]})$\footnote{We use a different convention from that of \cite{Sta23}, namely, $L_1=2L_1^{\text{Sta}}$, to give a clear comparison with our setting.}, 
which partly captures a non-orientable analogue of Mirzakhani's recursion. His approach is closely connected to matrix models and is also related, to some extent, to refined topological recursion.

Although Stanford's approach exhibits several interesting aspects, a crucial missing feature is that $\mathcal{V}_{g,n}^{\epsilon}(L_{[n]})\neq {V}_{g,n}^{\epsilon}(L_{[n]})$ as stated in \cite{Sta23}. More concretely, their difference is expected to be of order $O(\epsilon)$. In addition, since $\mathcal{V}_{g,n}^{\epsilon}(L_{[n]})$ are \emph{defined} by the proposed recursion formula, not geometrically in terms of the Norbury form and Gendulphe's regularisation, it is not guaranteed that $\mathcal{V}_{g,n}^{\epsilon}(L_{[n]})$ is symmetric under the permutation of elements in $L_{[n]}$ --- this issue is similar to \hyperref[item:RTR1]{RTR1} in Conjecture~\ref{conj:RTR}.

An important difference between Stanford's approach and refined topological recursion is that in \cite{Sta23}, unlike $\omega_{\frac12,1}$ in the refined spectral curve $\mathcal{S}^\mathfrak{b}_{\text{total}}$, a conventional volume $\mathcal{V}_{\frac12,1}^{\epsilon}(L_1)$ for M\"{o}bius strips is introduced as a part of the initial data as follows:
\[\mathcal{V}_{\frac12,1}^{\epsilon}(L_1)\coloneqq\frac{\theta(L_1>2\epsilon)}{2L_1\tanh\frac{L_1}{4}}.\]
Intuitively, from the perspective of hyperbolic geometry, the form of $\mathcal{V}_{\frac12,1}^{\epsilon}(L_1)$ suggests that the systole of 2-sided \emph{non-primitive} geodesics is bounded below by $\epsilon$, instead of the systole of 1-sided geodesics as in Gendulphe's regularisation. We note that by naively applying the dictionary \eqref{V and omega relation} for $\omega_{\frac12,1}$, we find:
\[V^{\epsilon,-}_{\frac12,1}(L_1)\coloneqq L_1^{-1}\hat{\mathcal{L}}_{L_1}^{-1}\omega_{\frac12,1}=\frac{\kappa(L_1,2\epsilon)}{2L_1\tanh\frac{L_1}{4}},\quad \text{ with }\;\; \kappa(L_1,\epsilon)=\begin{cases}
	1,\quad \text{ if } L>2\epsilon,\\
	\frac34,\quad \text{ if } L=2\epsilon,\\
	\frac12, \quad \text{ if } L<2\epsilon.
\end{cases}\]
Conversely, the Laplace transform of $L_1{V}_{\frac12,1}^{\epsilon,-}(L_1)$ diverges due to the singularity at $L_1=0$, but the Laplace transform of $L_1\mathcal{V}_{\frac12,1}^\epsilon(L_1)$ returns a finite differential which we denote by $\varpi_{\frac12,1}$:
\[\varpi_{\frac12,1}(z_1)\coloneqq\mathcal{L}_{z_1}L_1\mathcal{V}_{\frac12,1}^\epsilon(L_1)
=\left(\frac{e^{-2\epsilon z_1}}{2z_1}+\sum_{k\geq1}\frac{e^{-2(z_1+\frac{k}{2})\epsilon}}{z_1+\frac{k}{2}}\right)dz_1,\]
where the equality holds for ${\rm Re}(z_1)>0$ and the sum converges with the assumption that $\epsilon>0$. Following Stanford's approach, one may attempt to take $\varpi_{\frac12,1}$ as part of a refined spectral curve. However, although $\varpi_{\frac12,1}$  satisfies \hyperref[item:RSC1]{RSC1} and \hyperref[item:RSC2]{RSC2}, it fails to satisfy \hyperref[item:RSC3]{RSC3}.

\section{Euler characteristic $-1$ from refined topological recursion}\label{sec:evidence}
In order to provide supporting evidence, we prove that Conjecture~\ref{conj:RTR} and Question~\ref{conj:main} hold true, nontrivially,  for $\omega_{\frac12,2}$ and $\omega_{1,1}$. We emphasize that, although refined topological recursion requires tedious residue computations, these computations are entirely local in nature. Consequently, unlike the computation of the volumes $V_{g,n}^{\epsilon}$, the recursive construction of $\omega_{g,n}$ presents no conceptual obstructions in arbitrary topology.

\subsection{Two-bordered real projective planes}

\begin{proposition}\label{prop:1/2,2}
Conjecture~\ref{conj:RTR} holds for $\omega_{\frac12,2}$. Furthermore, it can be explicitly calculated as
\begin{align}
\omega_{\frac12,2}(z_1,z_2)=&2\mathfrak{b}\log\left(2\sinh\frac{\epsilon}{2}\right)\frac{dz_1dz_2}{z_1^2z_2^2}-\mathfrak{b}d_{z_1}\frac{\Delta_{z_1}\omega_{0,2}(z_1,z_2)}{4\omega_{0,1}(z_1)}-\mathfrak{b}d_{z_2}\frac{\Delta_{z_2}\omega_{0,2}(z_2,z_1)}{4\omega_{0,1}(z_2)}\nonumber\\
&+\frac{\mathfrak{b}}{2}\sum_{k\geq1}\frac{(-1)^k}{k}\left(\frac{dz_1}{(z_1-\frac{k}{2})^2}+\frac{dz_1}{(z_1+\frac{k}{2})^2}\right)\left(\frac{dz_2}{(z_2-\frac{k}{2})^2}+\frac{dz_2}{(z_2+\frac{k}{2})^2}\right)\label{explicit w_{1/2,2}}.
\end{align}
\end{proposition}
\begin{proof}
Note that if \eqref{explicit w_{1/2,2}} holds, then one can easily verify that Conjecture \ref{conj:RTR} holds. In particular, one can check that there is no pole as $z_1\to z_2$ and $z_1\to\frac{k}{2}$ for each $k\geq1$. Thus, we will derive such an explicit expression.

First of all, one has
\begin{equation*}
\left(\sum_{k\geq1}\Res_{z=\frac{k}{2}}+\sum_{i=1}^2\Res_{z=z_i}-\Res_{z=0}-\sum_{k\geq1}\Res_{z=-\frac{k}{2}}-\sum_{i=1}^2\Res_{z=-z_i}\right)\frac{\eta^z(z_1)}{4\omega_{0,1}(z)}\,\mathfrak{b}\,dx(z)\,d_z\frac{dx(z_2)}{(x(z)-x(z_2))^2}=0.
\end{equation*}
This is because the integrand is invariant under $\sigma:z\mapsto-z$, which implies that the residues at $z=\frac{k}{2}$ and $z=z_i$ cancel the residues at $z=-\frac{k}{2}$ and $z=-z_i$, respectively. Furthermore, the integrand is regular at $z=0$, so that there is no contribution from the residue at $z=0$. Then, by integrating by parts and also by utilising \eqref{Delta_0,2}, the formula for $\omega_{\frac12,2}$ can be brought into the following form:
\begin{multline}
\omega_{\frac12,2}(z_1,z_2)=\left(\sum_{k\geq1}\Res_{z=\frac{k}{2}}+\sum_{i=1}^2\Res_{z=z_i}-\Res_{z=0}-\sum_{k\geq1}\Res_{z=-\frac{k}{2}}-\sum_{i=1}^2\Res_{z=-z_i}\right)\\
\left(\frac{\eta^z(z_1)}{8\omega_{0,1}(z)}\Delta_z\omega_{\frac12,1}(z)\Delta_z\omega_{0,2}(z,z_2)-\mathfrak{b}\frac{\Delta_z\omega_{0,2}(z,z_1)\Delta_z\omega_{0,2}(z,z_2)}{8\omega_{0,1}(z)}\right)\label{1/2,2-1}.
\end{multline}
We denote the contribution from the first term in the second line of \eqref{1/2,2-1} by $\omega_{\frac12,2}^{(0)}$. The integrand of $\omega_{\frac12,2}^{(0)}$ is invariant under the involution $\sigma:z\mapsto-z$, hence all residues cancel except the one at $z=0$. Thus, we have
\begin{align*}
\omega_{\frac12,2}^{(0)}(z_1,z_2)=\mathfrak{b}\left(\epsilon-2\sum_{k\geq1}\frac{e^{-k\epsilon}}{k}\right)\frac{dz_1dz_2}{z_1^2z_2^2}=2\mathfrak{b}\log\left(2\sinh\frac{\epsilon}{2}\right)\frac{dz_1dz_2}{z_1^2z_2^2}.
\end{align*}

We next consider the contributions from the second term in the second line of \eqref{1/2,2-1}, which is independent of $\epsilon$. Since the integrand is anti-invariant under the involution $\sigma$, the residue at $z=0$ vanishes. We denote by $\omega_{\frac12,2}^{(1)}$ and $\omega_{\frac12,2}^{(2)}$ the sums of residues at $z=\pm z_1,\pm z_2$ and at $z=\pm\frac{k}{2}$ for $k\geq1$, respectively. Explicitly, one finds by tedious computations that
\begin{align*}
   \omega_{\frac12,2}^{(1)}(z_1,z_2)=&-\mathfrak{b}d_{z_1}\frac{\Delta_{z_1}\omega_{0,2}(z_1,z_2)}{4\omega_{0,1}(z_1)}-\mathfrak{b}d_{z_2}\frac{\Delta_{z_2}\omega_{0,2}(z_2,z_1)}{4\omega_{0,1}(z_2)},\\
     \omega_{\frac12,2}^{(2)}(z_1,z_2)=&\frac{\mathfrak{b}}{2}\sum_{k\geq1}\frac{(-1)^k}{k}\left(\frac{dz_1}{(z_1-\frac{k}{2})^2}+\frac{dz_1}{(z_1+\frac{k}{2})^2}\right)\left(\frac{dz_2}{(z_2-\frac{k}{2})^2}+\frac{dz_2}{(z_2+\frac{k}{2})^2}\right).
\end{align*}
This completes the proof.
\end{proof}


\begin{theorem}\label{thm:1/2,2} Question~\ref{conj:main} is answered affirmatively for $(g,n)=\big(\frac12,2\big)$. Explicitly,
 \begin{equation}
     L_1L_2\cdot V_{\frac12,2}^\epsilon(L_{[2]})= (1+b)^{\frac12}(\hat{\mathcal{L}}_{L_1}^{-1} \hat{\mathcal{L}}_{L_2}^{-1}).\; \omega_{\frac12,2}\Big|_{b=1}.
    \end{equation}
\end{theorem}
\begin{proof}
Let us assume $L_1\geq L_2$, and first consider $\mathcal{L}^{-1}_{L_2}.\;\omega_{\frac12,2}(z_1,z_2)$. Note that in general we have
\[\Res_{z=-u}\frac{1}{(z+u)^m}e^{z L} = \frac{1}{(m-1)!} L^{m-1} e^{-u L}.\]
Thus, thanks to Proposition \ref{prop:1/2,2}, the termwise inverse Laplace transform reduces to taking residues at the poles of $\omega_{\frac12,2}$. This can be explicitly computed as follows:
    \begin{align}
        \Res_{z_2=0}\omega_{\frac12,2}(z_1,z_2)e^{z_2L_2}=&2L_2\mathfrak{b}\log\left(2\sinh\frac{\epsilon}{2}\right)\frac{dz_1}{z_1^2}-\frac{\mathfrak{b}L_2^2}{2}\frac{dz_1}{z_1^2},\label{inverse of 1/2,2-1}\\
        \Res_{z_2=-\frac{k}{2}}\omega_{\frac12,2}(z_1,z_2)e^{z_2L_2}=&\mathfrak{b}L_2\frac{(-1)^k}{k}e^{-\frac{k}{2}L_2}\left(\frac{dz_1}{(z_1-\frac{k}{2})^2}+\frac{dz_1}{(z_1+\frac{k}{2})^2}\right),\label{inverse of 1/2,2-2}\\
        \Res_{z_2=-z_1}\omega_{\frac12,2}(z_1,z_2)e^{z_2L_2}=&-\mathfrak{b}\pi e^{-L_2 z_1}\left(\frac{L_2^2\,dz_1}{z_1\sin2\pi z_1}+\frac{L_2\,dz_1}{z_1^2\sin2\pi z_1}+\frac{2\pi L_2\cos2\pi z_1\,dz_1}{z_1\sin^22\pi z_1}\right).\label{inverse of 1/2,2-3}
    \end{align}
Note that the pole at $z=\frac{k}{2}$ in \eqref{inverse of 1/2,2-2} cancels the corresponding pole in \eqref{inverse of 1/2,2-3}; hence, $\mathcal{L}^{-1}_{L_2}.\;\omega_{\frac12,2}(z_1,z_2)$ has poles only at $z_1=0$ or $z_1=-\frac{k}{2}$.

Since we assumed $L_1\geq L_2$, one can now choose the contour for $\mathcal{L}^{-1}$ such that it encloses the origin and the negative real axis, regardless of the terms. From \eqref{inverse of 1/2,2-1} and \eqref{inverse of 1/2,2-2}, we have
\begin{align}
\Res_{z_1=0} \,e^{z_1L_1}\Res_{z_2=0}\omega_{\frac12,2}(z_1,z_2)\,e^{z_2L_2}=&2L_1L_2\mathfrak{b}\log\left(2\sinh\frac{\epsilon}{2}\right)-\frac{\mathfrak{b}L_1L_2^2}{2},\label{inverse of 1/2,2-4}\\
      \Res_{z_1=-\frac{k}{2}} \,e^{z_1L_1} \Res_{z_2=-\frac{k}{2}}\omega_{\frac12,2}(z_1,z_2)\,e^{z_2L_2}=&\mathfrak{b}L_1L_2\frac{(-1)^k}{k}e^{-\frac{k}{2}(L_1+L_2)}.\label{inverse of 1/2,2-5}
\end{align}
And for \eqref{inverse of 1/2,2-3}, we have contributions from $z_1=0$ and $z_1=-\frac{k}{2}$
    \begin{align*}
        \Res_{z_1=0}\,e^{z_1L_1} \Res_{z_2=-z_1}\omega_{\frac12,2}(z_1,z_2)e^{z_2L_2}=&-\frac{\mathfrak{b}}{2}L_1L_2(L_1-L_2),\\
        \Res_{z_1=-\frac{k}{2}}\,e^{z_1L_1} \Res_{z_2=-z_1}\omega_{\frac12,2}(z_1,z_2)e^{z_2L_2}=&\,\mathfrak{b}L_1L_2\frac{(-1)^k}{k}e^{-(L_1-L_2)\frac{k}{2}}.
    \end{align*}

Combining all these contributions, we obtain:
\begin{align*}
    (\hat{\mathcal{L}}_{L_1}^{-1} \hat{\mathcal{L}}_{L_2}^{-1}).\;\omega_{\frac12,2}=&\mathfrak{b}L_1L_2\left(2\log\left(2\sinh\frac{\epsilon}{2}\right)-\frac{L_1}{2}+\sum_{k\geq1}\frac{(-1)^k}{k}\left(e^{-(L_1+L_2)\frac{k}{2}}+e^{-(L_1-L_2)\frac{k}{2}}\right)\right)\nonumber\\
    =&\mathfrak{b}L_1L_2\left(2\log\left(2\sinh\frac{\epsilon}{2}\right)-\frac{L_1}{2}+\log\left(1+e^{-\frac{L_1+L_2}{2}}\right)-\log\left(1+e^{-\frac{L_1-L_2}{2}}\right)\right)\nonumber\\
    =&-\mathfrak{b}L_1L_2\left(\log\left(2\cosh\frac{L_1}{2}+2\cosh\frac{L_2}{2}\right)-2\log\left(2\sinh\frac{\epsilon}{2}\right)\right),
\end{align*}
where in the second equality we resummed the series in $k$ under the assumption that $L_1\geq L_2$. Thus, we arrive at the theorem by setting $b=1$, which means $\mathfrak{b}=-\frac{1}{\sqrt{2}}$. If we assume $L_1\leq L_2$, then we still obtain the same result by choosing a contour differently for the inverse Laplace transform of \eqref{inverse of 1/2,2-3}. Equivalently, one may repeat exactly the same argument by reversing the order of the inverse Laplace transforms, i.e.~applying $\mathcal{L}^{-1}_{L_1}\mathcal{L}^{-1}_{L_2}$, with an appropriate choice of contours. This explains how poles along the anti-diagonal should be treated.
\end{proof}

\subsection{One-bordered Klein bottles}\label{sec:RTR Klein bottle}
We will obtain the expression for $\omega_{1,1}$ in an explicit form. Since the derivation is heavily computational, let us first show the following lemma, which will be useful shortly:

\begin{lemma}\label{lem:w_{1,1}}
    \begin{align}
      &\mathfrak{b}\sum_{k\geq1}\Res_{z=\frac{k}{2}}\eta^z(z_1)d_{z_1}\frac{\Delta_{z}\omega_{\frac12,1}(z)}{4\omega_{0,1}(z)}\nonumber\\
      &=\mathfrak{b}^2\sum_{k\geq1}\frac{(-1)^k}{k}\left(\frac{dz_1}{(z_1-\frac{k}{2})^3}-\frac{dz_1}{(z_1+\frac{k}{2})^3}\right)+\mathfrak{b}^2\sum_{k\geq1}C_k^\epsilon\left(\frac{dz_1}{(z_1-\frac{k}{2})^2}+\frac{dz_1}{(z_1+\frac{k}{2})^2}\right),\nonumber
\end{align}
where
\begin{align}
C^\epsilon_k\coloneqq&2\frac{(-1)^k}{k}\log\left(2\sinh\frac{\epsilon}{2}\right)-\frac{(-1)^k}{k^2}+\frac{(-1)^k}{2k^2}(e^{k\epsilon}+e^{-k\epsilon})+\frac{(-1)^k}{k}\sum_{m=1}^{k-1}\frac{e^{m\epsilon}+e^{-m\epsilon}}{m}.\label{C for 1,1}
\end{align}
\end{lemma}
\begin{proof}
    Notice that as shown in e.g.~\cite[Section 2]{Osu23-1}, for any $a\in\mathbb{C}$ and $k\geq1$, we have
    \begin{equation}
        \Res_{z=a}\eta^z(z_1)\frac{dz}{(z-a)^{k+1}}=\frac{dz_1}{(z_1-a)^{k+1}}-(-1)^{k}\frac{dz_1}{(z_1+a)^{k+1}}.\label{eta-projection}
    \end{equation}
    Thus, the residue at $z=\frac{k}{2}$ is determined by expanding $d_z\frac{\Delta_{z}\omega_{\frac12,1}(z)}{4\omega_{0,1}(z_1)}$ at $z=\frac{k}{2}$. The expansion can be computed explicitly as follows
    \begin{align}
        d_z\frac{\Delta_{z}\omega_{\frac12,1}(z)}{4\omega_{0,1}(z_1)}=\mathfrak{b}^2\frac{(-1)^k}{k}\frac{dz}{\left(z-\frac{k}{2}\right)^3}+\mathfrak{b}^2\frac{\tilde C_k^\epsilon}{\left(z-\frac{k}{2}\right)^2}+\mathcal{O}(1),\label{expansion of delta}
    \end{align}
    where
     \begin{align}
       \tilde C_k^\epsilon =&\frac{(-1)^k}{k}\epsilon-\frac{(-1)^k}{2k^2}-\frac{(-1)^k}{2k^2}e^{-2\epsilon k}+\frac{1}{k^2}+\frac{(-1)^k}{2}\frac{e^{\epsilon k}-e^{-\epsilon k}}{k^2}\nonumber\\
&-2\sum_{i\geq1}^{i\neq k}\frac{(-1)^i}{(i-k)(k+i)}
-(-1)^k\sum_{j\geq1}^{j\neq k}\left(\frac{e^{-(j+k)\epsilon}}{k(j+k)}+\frac{e^{-(j-k)\epsilon}}{k(j-k)}\right).\label{tilde C_k^epsilon}
    \end{align}

In order to justify the expansion, recall that
    \begin{align*}
         \frac{1}{4\omega_{0,1}(z)}=&-\frac{1}{4zdz}\sum_{i\in\mathbb{Z}}\frac{(-1)^i}{z-\frac{i}{2}},\\
          \Delta_z\omega_{\frac12,1}(z)=&\mathfrak{b}\sum_{j\geq0}\frac{1}{1+\delta_{j0}}\frac{e^{2(z-\frac{j}{2})\epsilon}}{z-\frac{j}{2}}dz-\mathfrak{b}\sum_{j\leq0}\frac{1}{1+\delta_{j0}}\frac{e^{-2(z-\frac{j}{2})\epsilon}}{z-\frac{j}{2}}dz,
    \end{align*}
        where we rearranged the summation indices for convenience. One has to expand the derivative of the product of the above two at $z=\frac{k}{2}$ for $k>0$. If we keep the indices $(i,j)$ of the summand to indicate contributions, the term at $(i,j)=(k,k)$ gives the third order pole in  \eqref{expansion of delta} and also $\frac{(-1)^k}{k}\epsilon$ in $\tilde C^\epsilon_k$. For the other terms in the first line of $\tilde C^\epsilon_k$, $-\frac{(-1)^k}{2k^2}$ comes from $(i,j)=(-k,k)$, $\frac{(-1)^k}{2k^2}e^{-2\epsilon k}$ from $(i,j)=(k,-k)$, $\frac{1}{k^2}$ from $(i,j)=(0,k)$, and $\frac{(-1)^k}{2}\frac{e^{\epsilon k}-e^{-\epsilon k}}{k^2}$ from $(i,j)=(k,0)$. In the second line, the first summation term comes from $(i,j)=(i,k)$, with $i\neq\pm k,0$, and finally the second summation term from $(i,j)=(k,j)$, with $j\neq \pm k,0$. 

One can further simplify $\tilde C_k^\epsilon$ by a combinatorial manipulation (Lemma~\ref{lem:combinatorics0}) to show that $\tilde C_k^\epsilon=C_k^\epsilon$, and this completes the proof. 
\end{proof}

With this lemma, one can derive the explicit expression for  $\omega_{1,1}$ as below: 

\begin{proposition}\label{prop:1,1}
Conjecture~\ref{conj:RTR} holds for $\omega_{1,1}$. Furthermore, it can be explicitly calculated as
\begin{align}
\omega_{1,1}(z_1)=&\frac{1+5\mathfrak{b}^2}{8}\left(\frac{dz_1}{z_1^4}+\frac{2\pi^2}{3}\frac{dz_1}{z_1^2}\right)-\mathfrak{b}^2\left(\frac{\pi^2}{3}-2\log^2\sinh\frac{\epsilon}{2}\right)\frac{dz_1}{z_1^2}
-\mathfrak{b}d_{z_1}\frac{\Delta_{z_1}\omega_{\frac12,1}(z_1)}{4\omega_{0,1}(z_1)}\nonumber\\
&+\mathfrak{b}^2\sum_{k\geq1}\frac{(-1)^k}{k}\left(\frac{dz_1}{(z_1-\frac{k}{2})^3}-\frac{dz_1}{(z_1+\frac{k}{2})^3}\right)+\mathfrak{b}^2\sum_{k\geq1}C_k^\epsilon\left(\frac{dz_1}{(z_1-\frac{k}{2})^2}+\frac{dz_1}{(z_1+\frac{k}{2})^2}\right)\label{explicit w_{1,1}}.
\end{align}
\end{proposition}

\begin{proof}
Let us denote by $\omega_{1,1}^{(0)}$ the residue at $z=0$ of the refined topological recursion formula. One has to simply expand the integrand at $z=0$, which reads: 
\begin{align*}
\omega_{1,1}^{(0)}(z_1)=&\frac{1+5\mathfrak{b}^2}{8}\left(\frac{dz_1}{z_1^4}+\frac{2\pi^2}{3}\frac{dz_1}{z_1^2}\right)\nonumber\\
&-\mathfrak{b}^2\left(-\frac{\epsilon^2}{2}+2\sum_{k\geq1}\left(\frac{1}{k^2}+\frac{\epsilon e^{-\epsilon k}}{k}\right)-2\left(\sum_{k\geq1}\frac{e^{-\epsilon k}}{k}\right)^2\right)\frac{dz_1}{z_2^2}\nonumber\\
=&\frac{1+5\mathfrak{b}^2}{8}\left(\frac{dz_1}{z_1^4}+\frac{2\pi^2}{3}\frac{dz_1}{z_1^2}\right)-\mathfrak{b}^2\left(\frac{\pi^2}{3}-\log^2\sinh\frac{\epsilon}{2}\right)\frac{dz_1}{z_2^2}.
\end{align*}
It is also not difficult to evaluate the residues at $z=\pm z_1$, whose contributions are denoted by $\omega_{1,1}^{(1)}$, and we have (c.f.~\cite[Lemma 3.11]{Osu23-1}):
\begin{equation}
\omega_{1,1}^{(1)}(z_1)=-\frac{{\rm Rec}_{1,1}(z_1)-{\rm Rec}_{1,1}(-z_1)}{2\omega_{0,1}(z_1)}=-\mathfrak{b}d_{z_1}\frac{\Delta_{z_1}\omega_{\frac12,1}(z_1)}{4\omega_{0,1}(z_1)}.\label{w_{1,1}-1}
\end{equation}

We denote the remaining contributions by $\omega_{1,1}^{(2)}$. Notice that the invariant part of the integrand vanishes since we take the difference between the residues at $z=\pm\frac{k}{2}$ for each $k\geq1$. Therefore, we can simplify it to
\begin{equation*}
\omega_{1,1}^{(2)}(z_1)=\sum_{k\geq1}\left(\Res_{z=\frac{k}{2}}-\Res_{z=-\frac{k}{2}}\right)\frac{\eta^z(z_1)}{4\omega_{0,1}(z)}{\rm Rec}_{1,1}(z)=\mathfrak{b}\sum_{k\geq1}\Res_{z=\frac{k}{2}}\eta^z(z_1)d_{z_1}\frac{\Delta_{z}\omega_{\frac12,1}(z)}{4\omega_{0,1}(z)},
\end{equation*}
where we used the second equality of \eqref{w_{1,1}-1} at the second equality of this equation. This is already computed in Lemma~\ref{lem:w_{1,1}}, and combining all $\omega_{1,1}^{(0)}$, $\omega_{1,1}^{(1)}$, $\omega_{1,1}^{(2)}$, one obtains the explicit expression of $\omega_{1,1}$, from which it is straightforward to check that $\omega_{1,1}$ has no residue and has poles only at the expected locations.
\end{proof}

We note that the last summation involving $C_k^\epsilon$ diverges, as $C_k^\epsilon$ contains an exponential factor $e^{k\epsilon}$, while the denominator decays only polynomially as $k^{-2}$. A crucial feature of the termwise inverse Laplace transform $\hat{\mathcal{L}}^{-1}$ is that it transforms $(z_1+\frac{k}{2})^{-2}$ into $L_1 \exp(-\frac{k}{2}L_1)$, which ensures convergence of the resulting series provided that $L_1>2\epsilon$. Furthermore, after applying $\hat{\mathcal{L}}^{-1}$ and resummation, one can analytically continue the result in the variable $L_1$ to all of $\mathbb{R}_+$, thereby removing the restriction $L_1>2\epsilon$. This phenomenon motivates the property \hyperref[item:RTR0]{RTR0} in Conjecture~\ref{conj:RTR}.

We now prove that $\omega_{1,1}$ and $V^\epsilon_{1,1}$ are related by the formula in Question~\ref{conj:main}.

\begin{theorem}\label{thm:1,1}
 \begin{equation}
     L_1\cdot V_{1,1}^\epsilon(L_1)= (1+b)\hat{\mathcal{L}}^{-1}_{L_1}.\;\omega_{1,1}\Big|_{b=1},
    \end{equation} 
    by analytic continuation on the right-hand side with respect to $L_1$ as described above.
\end{theorem}
\begin{proof}
Let us assume $L_1>2\epsilon$. Thanks to Proposition \ref{prop:1,1}, it is easy to check that
    \begin{align}
        \Res_{z=0}\omega_{1,1}(z)\,e^{z L_1}=&\frac{1+5\mathfrak{b}^2}{48}L_1^3-\mathfrak{b}^2L_1^2\log\left(2\sinh\frac{\epsilon}{2}\right)+\frac{\pi^2}{12}L_1+\mathfrak{b}^2L_1\left(\frac{\pi^2}{12}+2\log^2\left(2\sinh\frac{\epsilon}{2}\right)\right),\label{vol at z=0}\\
        \Res_{z=-\frac{k}{2}}\omega_{1,1}(z)\,e^{z L_1}=&\mathfrak{b}^2\left(-\frac{(-1)^k}{k}L_1^2+2L_1C_k^\epsilon \right)e^{-\frac{k}{2}L}.\label{vol at z=-k/2}
    \end{align}
Note that the second term in \eqref{vol at z=0} comes from the expansion of $\omega_{1,1}^{(1)}$ (the last term in the first line of \eqref{explicit w_{1,1}}) at $z=0$, which reads:
    \begin{equation*}
     -  \mathfrak{b}d_{z}\frac{\Delta_{z}\omega_{\frac12,1}(z)}{4\omega_{0,1}(z)}=-\mathfrak{b}^2 \left(\epsilon-\sum_{k\geq1}\frac{2}{k}e^{-k\epsilon}\right)\frac{dz}{z^3}+\mathcal{O}(1)=-2\mathfrak{b}^2\log\left(2\sinh\frac{\epsilon}{2}\right)\frac{dz}{z^3}+\mathcal{O}(1),
    \end{equation*}
where the second equality holds due to the assumption that $\epsilon\in\mathbb{R}_{>0}$. Also, the residue of $\omega_{1,1}^{(1)}$ at $z=\frac{k}{2}$ cancels that of $\omega_{1,1}^{(2)}$, whereas the residue at $z=-\frac{k}{2}$ doubles. This is why we have $2L_1C_k^\epsilon$, not $L_1C_k^\epsilon$, in the second term in \eqref{vol at z=-k/2}.

We can show that $U_k^\epsilon$ in \eqref{U for 1,1} and $C_k^\epsilon$ in \eqref{C for 1,1} satisfy $(-1)^kU_k^\epsilon=2C_k^\epsilon$ by nontrivial combinatorial manipulations (Lemma~\ref{lem:combinatorics2}). Thus we have:
\begin{align*}
\hat{\mathcal{L}}_{L_1}^{-1}.\;\omega_{1,1}(z_1)=&\frac{1+5\mathfrak{b}^2}{48}L_1^3-\mathfrak{b}^2L_1^2\log\left(2\sinh\frac{\epsilon}{2}\right)+\frac{\pi^2}{12}L_1+\mathfrak{b}^2L_1\left(\frac{\pi^2}{12}+2\log^2\left(2\sinh\frac{\epsilon}{2}\right)\right)\nonumber\\
&-\mathfrak{b}^2L_1^2\sum_{k\geq1}\frac{(-1)^k}{k}e^{-\frac{k}{2}L_1}+\mathfrak{b}^2L_1\sum_{k\geq1}(-1)^k U^\epsilon_ke^{-\frac{k}{2}L_1}.
\end{align*}
Then, it can be easily seen that the theorem holds when $b=1$ ($\mathfrak{b}^2=\frac12$) by comparing it with Proposition \ref{prop:expansion of V_KB} --- do not forget to add $V_{1,1}^+(L_1)$ of \eqref{V_{0,3} and V_{1,1}}. Note that strictly speaking the equality holds when $L_1>2\epsilon$, where the expansion of $V_{1,1}^\epsilon$ is valid, as described in Proposition \ref{prop:expansion of V_KB}. Once we resum the expansion into a dilogarithm, then one can analytically continue with respect to $L_1$ to the entire $\mathbb{R}_+$, not only $L_1>2\epsilon$.
\end{proof}


\appendix

\section{Additional technical results and discussions}


\subsection{Combinatorial formulae for the $(1,1)$ case}\label{sec:formulae}

We present a collection of combinatorial lemmas that simplify computations arising in refined topological recursion, in particular those needed to establish the correspondence with the geometric side in the $(g,n)=(1,1)$ case. In practice, we need to substitute $X = e^{\frac{\epsilon}{2}}$ to use the lemmas below in Section~\ref{sec:RTR Klein bottle}.
The first result is used in Lemma~\ref{lem:w_{1,1}}.
\begin{lemma}\label{lem:combinatorics0}
For $X>1$,
\begin{align*} 
    &2\sum_{i\geq1}^{i\neq k}\frac{(-1)^i}{(i-k)(i+k)}
+(-1)^k\sum_{j\geq1}^{j\neq k}\left(\frac{X^{-2(j+k)}}{k(j+k)}+\frac{X^{-2(j-k)}}{k(j-k)}\right)\nonumber\\
    &=\frac{1}{k^2}-\frac{(-1)^k}{k}\left(\frac{X^{-4k}-1}{2k}+2\log\left(1-\frac{1}{X^{2}}\right)+\frac{1}{kX^{2k}}+\sum_{j=1}^{k-1}\frac{1}{j}\left(X^{2j}+\frac{1}{X^{2j}}\right)\right).
\end{align*}
\end{lemma}

\begin{proof}
First, notice that
\begin{equation}
    \sum_{i\geq1}^{i\neq k}\frac{(-1)^i}{(i-k)}=(-1)^k\left(\sum_{i\geq1}\frac{(-1)^i}{i}-\sum_{i=1}^{k-1}\frac{(-1)^i}{i}\right).\label{comb0-1}
\end{equation}
Also,
\begin{align*}
     \sum_{i\geq1}^{i\neq k}\frac{(-1)^i}{(i+k)}=&-\frac{(-1)^k}{2k}+\sum_{i\geq1}\frac{(-1)^i}{(i+k)}=-\frac{(-1)^k}{2k}+(-1)^k\left(\sum_{i\geq1}\frac{(-1)^i}{i}-\sum_{i=1}^{k}\frac{(-1)^i}{i}\right),
\end{align*}
where the finite sum stops at $i=k$ in contrast to the finite sum in \eqref{comb0-1}. Thus, we have
\begin{align*}
   2\sum_{i\geq1}^{i\neq k}\frac{(-1)^i}{(i-k)(i+k)}=&\frac{1}{k}\sum_{i\geq1}^{i\neq k}\left(\frac{(-1)^i}{(i-k)}-\frac{(-1)^i}{(i+k)}\right)=\frac{(-1)^k}{2k^2}+\frac{1}{k^2}.
\end{align*}

Next, one can manipulate the following sum as below:
\begin{equation}
    \sum_{j\geq1}^{j\neq k}\frac{X^{-2(j-k)}}{(j-k)}=\left(\sum_{j\geq1}\frac{X^{-2j}}{j}-\sum_{j=1}^{k-1}\frac{X^{2j}}{j}\right)=-\log\left(1-\frac{1}{X^2}\right)-\sum_{j=1}^{k-1}\frac{X^{2j}}{j},\label{comb0-2}
\end{equation}
where the second equality holds because the first sum is a convergent series when $X>1$. Similarly, we have
\begin{align*}
    \sum_{j\geq1}^{j\neq k}\frac{X^{-2(j+k)}}{(j+k)}=-\frac{1}{2k}\frac{1}{X^{4k}}+  \sum_{j\geq1}\frac{X^{-2(j+k)}}{(j+k)}=-\frac{1}{2k}\frac{1}{X^{4k}}-\log\left(1-\frac{1}{X^2}\right)-\sum_{j=1}^{k}\frac{1}{j}\frac{1}{X^{2j}},
\end{align*}
where the finite sum stops at $j=k$ in contrast to the finite sum in \eqref{comb0-2}. Combining all of the above, one arrives at the lemma.
\end{proof}

The following lemma is necessary to show the forthcoming lemma:
\begin{lemma}\label{lem:combinatorics1}
    \begin{equation}
        \frac{1}{k}\left(X^{2k}+\frac{1}{X^{2k}}\right)-\frac{2}{k}=\sum_{j=1}^k\frac{1}{j}\binom{k+j-1}{2j-1}\left(X-\frac{1}{X}\right)^{2j}.
    \end{equation}
\end{lemma}
\begin{proof}
    Let us consider $\phi_k(X)\coloneqq\frac{1}{k}\left(X^{2k}+\frac{1}{X^{2k}}\right)$. Then there exists a set of constants $c_j$ such that
    \begin{equation*}
        \phi_k(X)=\sum_{j=0}^kc_j\left(X-\frac{1}{X}\right)^{2j}.
    \end{equation*}
    It is clear that $\phi_k(X)$ satisfies $\left(X\frac{d}{dX}\right)^2\phi_k(X)=4k^2\phi_k(X)$, and this differential equation induces the following recursive equation for $c_j$ starting from $c_k$:
    \begin{equation*}
      \forall\;\; 0\leq j\leq k-1,\qquad  (2j+1)(2j+2)c_{j+1}=(k-j)(k+j)c_j.
    \end{equation*}
    This can be solved by induction on $k-j$ and one obtains:
    \begin{equation*}
        c_j=\frac{k}{j}\binom{k+j-1}{2j-1}c_k.
    \end{equation*}
    We note that the above equation holds even when $j=0$ which gives $c_0=2c_k$. Finally, the boundary condition $\phi(\pm1)=\frac{2}{k}$ fixes $c_k=\frac{1}{k}$.
\end{proof}
The lemma below is used in Theorem \ref{thm:1,1}.
\begin{lemma}\label{lem:combinatorics2}
    \begin{equation*}
         \sum_{j=1}^{k-1}\frac{4}{j}+k\sum_{j=1}^k\frac{1}{j^2}\binom{k+j-1}{2j-1}\left(X-\frac{1}{X}\right)^{2j}=-\frac{2}{k}+\frac{1}{k}\left(X^{2k}+\frac{1}{X^{2k}}\right)+2\sum_{j=1}^{k-1}\frac{1}{j}\left(X^{2j}+\frac{1}{X^{2j}}\right).
    \end{equation*}
\end{lemma}
\begin{proof}
    We proceed by induction on $k$. When $k=1$, it is straightforward to check that the lemma holds. Let $k\geq1$ and we assume that the lemma holds for all $k'\leq k$. Then, by utilising Lemma~\ref{lem:combinatorics1}, one finds:
    \begin{align*}
        &-\frac{2}{k+1}+\frac{1}{k+1}\left(X^{2k+2}+\frac{1}{X^{2k+2}}\right)+2\sum_{j=1}^{k}\frac{1}{j}\left(X^{2j}+\frac{1}{X^{2j}}\right)\nonumber\\
        &= -\frac{2}{k+1}+\frac{1}{k+1}\left(X^{2k+2}+\frac{1}{X^{2k+2}}\right)+\frac{1}{k}\left(X^{2k}+\frac{1}{X^{2k}}\right)\nonumber\\
        &\;\;\;\;+\frac{2}{k}+\sum_{j=1}^{k-1}\frac{4}{j}+k\sum_{j=1}^k\frac{1}{j^2}\binom{k+j-1}{2j-1}\left(X-\frac{1}{X}\right)^{2j}\nonumber\\
        &=\sum_{j=1}^{k+1}\frac{1}{j}\binom{k+j}{2j-1}\left(X-\frac{1}{X}\right)^{2j}+\frac{2}{k}+\sum_{j=1}^k\frac{1}{j}\binom{k+j-1}{2j-1}\left(X-\frac{1}{X}\right)^{2j}\nonumber\\
        &\;\;\;\;+\frac{2}{k}+\sum_{j=1}^{k-1}\frac{4}{j}+k\sum_{j=1}^k\frac{1}{j^2}\binom{k+j-1}{2j-1}\left(X-\frac{1}{X}\right)^{2j}\nonumber\\
        &=\sum_{j=1}^{k}\frac{4}{j}+\frac{1}{k+1}\left(X-\frac{1}{X}\right)^{2k+2}\nonumber\\
        &\;\;\;\;+\sum_{j=1}^k\frac{1}{j^2}\left(j\binom{k+j}{2j-1}+j\binom{k+j-1}{2j-1}+k\binom{k+j}{2j-1}\right)\left(X-\frac{1}{X}\right)^{2j},
    \end{align*}
    where the first equality is due to the induction hypothesis, and we used Lemma~\ref{lem:combinatorics1} at the second equality. Finally, expressing the binomial coefficients in terms of factorials simplifies the last line and completes the proof.
\end{proof}

\subsection{Limitations of the Laplace transform approach}\label{sec:subtleties_in_question_ref_conj_main}

The relation between Mirzakhani's recursion and the Eynard--Orantin topological recursion was originally proven in \cite{EO07-2} by direct computations with several clever manipulations. One can also prove it in terms of Airy structures and Virasoro constraints \cite{KS17}. Since the latter is not known in the refined setting at the moment of writing, a realistic approach for Question~\ref{conj:main} is to extend the original proof by Eynard and Orantin to our problem.

Let us first briefly explain how the relation between the topological recursion correlators $\omega_{g,n}$ and the Weil--Petersson volumes $V_{g,n}$ was originally discovered in the unrefined/oriented setting \cite{EO07-2}. 
For $(g,n)$ with $2g-2+n>0$, we symbolically denote the Eynard--Orantin recursive formula by
\begin{equation}
    \omega_{g,n}=\text{EO}\Big(\omega_{0,1},\omega_{0,2};\{\omega_{h,m}\}_{0<2h-2+m<2g-2+n}\Big).\label{symbolicEO}
\end{equation}
Let us next define $U_{g,n}(L_{[n]})$ by
\[U_{g,n}(L_{[n]})\coloneqq \prod_{i=1}^nL_i^{-1}\cdot(\mathcal{L}_{L_1}^{-1}\cdots \mathcal{L}_{L_n}^{-1}).\;\omega_{g,n},\]
where $\mathcal{L}^{-1}$ is the inverse Laplace transform, not the termwise one $\hat{\mathcal{L}}^{-1}$. Since the $\omega_{g,n}$ have poles only at ramification points, it can be easily shown that $U_{g,n}(L_{[n]})$ are polynomials in $L_1,...,L_n$, and hence, taking the inverse is straightforward:
\[\omega_{g,n}=(\mathcal{L}_{z_1}\cdots\mathcal{L}_{z_n}).\;\prod_{i=1}^n L_i\cdot U_{g,n}(L_{[n]}).\]

One can check by explicit computations that $U_{0,3}=V_{0,3}$ and $U_{1,1}=V_{1,1}$. For $2g-2+n>1$, the symbolic formula \eqref{symbolicEO} can be written as
\begin{equation}
  U_{g,n}(L_{[n]})=  \prod_{i=1}^nL_i^{-1}\cdot(\mathcal{L}^{-1}_{L_i}\cdots\mathcal{L}^{-1}_{L_n}).\;\text{EO}\Big(\omega_{0,1},\omega_{0,2};\{(\mathcal{L}_{z_1}\cdots\mathcal{L}_{z_m}).\prod_{i=1}^mL_i\;U_{h,m}\}_{0<2h-2+m<2g-2+n}\Big).\label{EO proof1}
\end{equation}
This provides a recursion formula for $U_{g,n}$, induced from the Eynard--Orantin recursion for $\omega_{g,n}$. It was then shown in \cite{EO07-2} that the recursion for $U_{g,n}$ precisely coincides with Mirzakhani's recursion, implying $U_{g,n}=V_{g,n}$ for all $(g,n)$. The nontrivial manipulations in \cite{EO07-2} transform the dependence of $\omega_{0,1}$ and $\omega_{0,2}$ in \eqref{EO proof1} to the functions $R$ and $D$ in Mirzakhani's recursion.

There are a few obstacles to extending this idea to the refined setting, even if one assumes the properties of Conjecture \ref{conj:RTR}. Besides technical difficulties, one main issue is that the manipulations in \cite{EO07-2} rely on the facts that, in the unrefined setting, $V_{g,n}$ and $\omega_{g,n}$ are in a one-to-one correspondence through the Laplace transform and its inverse. In the non-orientable setting, however, it turns out that $\omega_{\frac12,1}$ is written in the following integral form, with the assumption that $\text{Re}(z)>0$,
\[
     \omega_{\frac12,1}(z)=\frac{b}{2}\left(-\int_{2\epsilon}^\infty \frac{1}{\tanh\frac{q}{4}}e^{-zq}dq-\int_{0}^{2\epsilon}\frac{\sinh zq}{\tanh\frac{q}{4}}dq\right)dz.
\]
Since this cannot be written as a Laplace transform (i.e.~the $z$-dependence appears only in $e^{-zq}$), one cannot directly apply some of the clever manipulations of \cite{EO07-2}. Furthermore, recall that $\hat{\mathcal{L}}_{L_1}\omega_{1,1}$ is divergent when $0<L\leq2\epsilon$, which means that $\omega_{1,1}$ cannot be inverted to $V_{1,1}^{\epsilon}(L_1)$. This suggests that a recursive structure on $\omega_{g,n}$ cannot be readily translated into one for $V_{g,n}^{\epsilon}(L_{[n]})$, consistent with the geometric difficulties discussed above.

\addcontentsline{toc}{section}{References}
\printbibliography
\end{document}